\documentclass[psamsfonts,reqno]{amsart}
\usepackage{amscd,amssymb,latexsym,graphics}
\usepackage{epstopdf}

\theoremstyle{plain}

\newtheorem{lemma}{Lemma}[section]
\newtheorem{theorem}[lemma]{Theorem}
\newtheorem{proposition}[lemma]{Proposition}
\newtheorem{corollary}[lemma]{Corollary}
\newtheorem{examplepf}[lemma]{Example}

\newtheorem*{sclaim}{Claim}

\newtheorem*{stat}{\name}
\newcommand{\name}{testing}

\theoremstyle{definition}
\newtheorem{definition}[lemma]{Definition}
\newtheorem{example}[lemma]{Example}

\theoremstyle{remark}
\newtheorem{remark}[lemma]{Remark}
\newtheorem{notation}[lemma]{Notation}

\newcommand{\qedc}{{\qed}~{\rm Claim~{\theclaim}.}}
\newcommand{\qedsc}{{\qed}~{\rm Claim.}}

\newenvironment{scproof}
{\begin{proof}[Proof of Claim.]}
{\qedsc\renewcommand{\qed}{}\end{proof}}

\numberwithin{equation}{section}
\numberwithin{figure}{section}

\newcommand{\cc}[1]{\hfill\kern.7em#1\kern.7em\hfill}

\newcommand{\bve}[2]{\Vert{#1}\nobreak=\nobreak{#2}\Vert}
\newcommand{\bvo}[2]
{[\![{#1}\nobreak\leqslant\nobreak{#2}]\!]}
\newcommand{\bveo}[2]
{[\![{#1}\nobreak=\nobreak{#2}]\!]}
\newcommand{\bvi}[2]{\Vert{#1}\nobreak\leqslant\nobreak{#2}\Vert}
\newcommand{\bvp}[2]{\bvi{{#1}}{{#2}}^+}
\newcommand{\bvm}[2]{\bvi{{#1}}{{#2}}^-}
\newcommand{\bvpm}[2]{\bvi{{#1}}{{#2}}^{\pm}}

\newcommand{\bvP}[2]{\bvi{{#1}}{{#2}}_{\bP}}
\newcommand{\bvQ}[2]{\bvi{{#1}}{{#2}}_{\bQ}}
\newcommand{\pup}[1]{\textup{(}{#1}\textup{)}}
\newcommand{\bl}{\bullet}

\newcommand{\into}{\hookrightarrow}

\newcommand{\inc}{\mathbin{\Vert}}
\newcommand{\Pow}{\mathfrak{P}}

\newcommand{\lerc}{\leq_{\mathrm{rc}}}
\newcommand{\leint}{\leq_{\mathrm{int}}}

\newcommand{\lecov}{\leq_{\mathrm{cov}}}
\newcommand{\ledb}{\leq_{\mathrm{db}}}

\DeclareMathOperator{\J}{J}

\newcommand{\hgt}{\mathrm{height}}

\newcommand{\cA}{\mathcal{A}}
\newcommand{\cB}{\mathcal{B}}
\newcommand{\cI}{\mathcal{I}}
\newcommand{\cC}{\mathcal{C}}
\newcommand{\cD}{\mathcal{D}}

\newcommand{\xA}{\mathbf{A}}
\newcommand{\xB}{\mathbf{B}}

\newcommand{\DSLat}{\mathbf{DSLat}}
\newcommand{\DSLatemb}{\mathbf{DSLat}^{\mathrm{emb}}}
\newcommand{\DLat}{\mathbf{DLat}}
\newcommand{\DLatemb}{\mathbf{DLat}^{\mathrm{emb}}}
\newcommand{\Lat}{\mathbf{Lat}}
\newcommand{\VPmeas}{\mathbf{VPMeas}}

\newcommand{\es}{\varnothing}
\newcommand{\jz}{$\langle\vee,0\rangle$}
\newcommand{\jzu}{$\langle\vee,0,1\rangle$}
\newcommand{\jzs}{\jz-semi\-lat\-tice}
\newcommand{\jzus}{\jzu-semi\-lat\-tice}
\newcommand{\jzh}{\jz-ho\-mo\-mor\-phism}
\newcommand{\jzuh}{\jzu-ho\-mo\-mor\-phism}
\newcommand{\jze}{\jz-em\-bed\-ding}
\newcommand{\jzue}{\jzu-em\-bed\-ding}
\newcommand{\mz}{$\langle\wedge,0\rangle$}

\newcommand{\mzs}{\mz-semi\-lat\-tice}

\newcommand{\js}{join-sem\-i\-lat\-tice}
\newcommand{\ms}{meet-sem\-i\-lat\-tice}
\newcommand{\jirr}{join-ir\-re\-duc\-i\-ble}

\newcommand{\contr}{a contradiction}

\newcommand{\res}{\mathbin{\restriction}}
\newcommand{\lc}[2]{{#1}_{({#2})}}
\newcommand{\uc}[2]{{#1}^{({#2})}}

\newcommand{\set}[1]{\{{#1}\}}

\newcommand{\setm}[2]{\set{#1\mid#2}}

\newcommand{\seq}[1]{\langle{#1}\rangle}
\newcommand{\seqm}[2]{\seq{{#1}\mid{#2}}}
\newcommand{\famm}[2]{({#1}\mid{#2})}
\newcommand{\Famm}[2]{\bigl({#1}\mid{#2}\bigr)}

\newcommand{\bP}{\boldsymbol{P}}
\newcommand{\bQ}{\boldsymbol{Q}}
\newcommand{\bR}{\boldsymbol{R}}

\newcommand{\id}{\mathrm{id}}

\newcommand{\ol}[1]{\,\overline{\!#1}}

\DeclareMathOperator{\Conc}{Con_c}

\DeclareMathOperator{\Idc}{Id_c}

\newcommand{\eps}{\varepsilon}

\newcommand{\ba}{\boldsymbol{a}}
\newcommand{\bb}{\boldsymbol{b}}
\newcommand{\bc}{\boldsymbol{c}}

\newcommand{\bx}{\boldsymbol{x}}
\newcommand{\by}{\boldsymbol{y}}

\newcommand{\sx}{\mathsf{x}}
\newcommand{\sy}{\mathsf{y}}

\begin{document}

\author[F.~Wehrung]{Friedrich Wehrung}
\address{CNRS, UMR 6139\\
D\'epartement de Math\'ematiques, BP 5186\\
Universit\'e de Caen, Campus 2\\
14032 Caen cedex\\
France}
\email{wehrung@math.unicaen.fr}
\urladdr{http://www.math.unicaen.fr/\~{}wehrung}
\date{\today}

\title{Poset representations of distributive
semilattices}
\date{\today}

\begin{abstract}
We prove that for every distributive \jzs\ $S$, there are a
\ms\ $P$ with zero and a map $\mu\colon P\times P\to S$ such
that $\mu(x,z)\leq\mu(x,y)\vee\mu(y,z)$ and $x\leq y$ implies
that $\mu(x,y)=0$, for all $x,y,z\in P$, together with the
following conditions:
\begin{enumerate}
\item[(P1)] $\mu(v,u)=0$ implies that $u=v$, for all $u\leq v$ in $P$.

\item[(P2)] For all $u\leq v$ in~$P$ and all $\ba,\bb\in S$, if
$\mu(v,u)\leq\ba\vee\bb$, then there are a positive integer $n$ and
a decomposition $u=x_0\leq x_1\leq\cdots\leq x_n=v$ such that either
$\mu(x_{i+1},x_i)\leq\ba$ or $\mu(x_{i+1},x_i)\leq\bb$, for each $i<n$.

\item[(P3)] The subset $\setm{\mu(x,0)}{x\in P}$ generates the
semilattice $S$.
\end{enumerate}
Furthermore, every finite, bounded subset of $P$ has a join, and
$P$ is bounded in case $S$ is bounded.
Furthermore, the construction is functorial on lattice-indexed
diagrams of finite distributive \jzus s.
\end{abstract}

\maketitle

\section{Introduction}\label{S:Intro}

\subsection{Origin of the problem}\label{Su:Origin}
The classical congruence lattice representation problem, usually
denoted by~CLP, asks whether every distributive \jzs\ is
isomorphic to the semilattice $\Conc L$ of all compact (i.e., finitely generated) congruences of some lattice $L$. (It is well-known, see \cite{FuNa42} or
\cite[Theorem~II.3.11]{GLT2}, that $\Conc L$ is a distributive
\jzs, for every lattice $L$.) This problem has finally been solved negatively by the author in~\cite{CLP}. This negative solution came out of a failed attempt
to extend to \emph{semilattices} the representation result of
distributive semilattices by \emph{posets} (i.e., partially ordered sets) stated in the Abstract. The purpose of the present paper is to give a proof of that result.

A first motivation for proving this result lies in its relation with congruence lattices of lattices, which we shall outline now.

\begin{definition}\label{D:PosetMeas}
Let $S$ be a \jzs\ and let $P$ be a poset. A map\linebreak
$\mu\colon P\times P\to S$ is a \emph{$S$-valued p-measure
on~$P$}, if $\mu(x,z)\leq\mu(x,y)\vee\mu(y,z)$ and $x\leq y$
implies $\mu(x,y)=0$, for all $x,y,z\in P$. The pair
$\seq{P,\mu}$ is a \emph{$S$-valued p-measured poset}.
\end{definition}

The inequality $\mu(x,z)\leq\mu(x,y)\vee\mu(y,z)$ will be
referred to as the \emph{triangular inequality}. The letter `p' in `p-measure' stands for `poset'.

\begin{notation}\label{Not:BVs}
We shall always denote by $P$, $Q$, $\dots$, the underlying
posets of p-measured posets $\bP$, $\bQ$, $\dots$\,. For a
p-measured poset $\bP=\seq{P,\mu}$, we shall often use the notation
$\bvP{x}{y}=\mu(x,y)$, for $x,y\in P$. Elements of the form
$\bvP{x}{y}$ will be called \emph{Boolean values}.
\end{notation}

A fundamental class of p-measures is given as follows. For a lattice~$L$, the map $\Theta^+\colon L\times L\to\Conc L$ defined by the rule
 \begin{equation}\label{Eq:Thetaplus}
 \Theta^+(x,y)=\text{least congruence of }L\text{ that identifies }x\vee y
 \text{ and }y,
 \end{equation}
for all $x,y\in L$, is obviously a $\Conc L$-valued p-measure on~$L$. Furthermore, it satisfies the following conditions:
\begin{enumerate}
\item $\Theta^+(v,u)=0$ implies that $u=v$, for all $u\leq v$ in $L$.

\item For all $u\leq v$ in~$L$ and all $\ba,\bb\in\Conc L$, if
$\Theta^+(v,u)\leq\ba\vee\bb$, then there are a positive integer $n$ and
a decomposition $u=x_0\leq x_1\leq\cdots\leq x_n=v$ such that either
$\Theta^+(x_{i+1},x_i)\leq\ba$ or $\Theta^+(x_{i+1},x_i)\leq\bb$, for each $i<n$.

\item The subset $\setm{\Theta^+(v,u)}{u\leq v\text{ in }L}$ generates the
semilattice $\Conc L$.
\end{enumerate}

Item (ii) above follows from the usual description of congruences in lattices, see, for example, \cite[Theorem~I.3.9]{GLT2}, while Items~(i) and~(iii) are trivial.

Hence, if a distributive \jzs\ $S$ is isomorphic to $\Conc L$ for some lattice~$L$, then there exists a $S$-valued p-measure $\mu\colon L\times L\to S$ satisfying (i)--(iii) above. Although we could prove in~\cite{CLP} that there may not exist such a lattice~$L$, the main result of the present paper is that the conclusion about p-measures persists. Conditions~(i) and~(ii) are the same as the conditions denoted by~(P1) and~(P2), respectively, in the Abstract, while~(P3) is a strengthening of~(iii).

Furthermore, unlike earlier representation results such as Gr\"atzer and Schmidt's representation theorem of algebraic lattices as congruence lattices of abstract algebras \cite{GrSc63}, our result \emph{characterizes} distributive algebraic lattices, see Corollary~\ref{C:PosetRepr} and Proposition~\ref{P:CondImplDistr}.

\subsection{Lifting objects and diagrams with respect to functors}
\label{Su:Lifting}

Most of the recent efforts at solving CLP have been aimed at lifting
not only individual (distributive) semilattices, but also \emph{diagrams} of
semilattices, with respect to the congruence semilattice functor
$\Conc$. They are based on the following lemma, proved by Ju.\,L. Ershov as the main theorem in Section~3 of the Introduction of his 1977 monograph~\cite{Ersh} and P. Pudl\'ak in his 1985 paper \cite[Fact~4, p.~Ê100]{Pudl}.

\begin{lemma}\label{L:ErsPudl}
Every distributive \jzs\ $S$ is the directed union of its finite distributive \jz-subsemilattices.
\end{lemma}

Because of this, lifting diagrams of distributive semilattices can be reduced to lifting diagrams of \emph{finite} distributive semilattices.

The formal definition of a lifting runs as follows. For categories $\cI$, $\cA$, $\cB$ and a functor $\Phi\colon\cA\to\cB$, a \emph{lifting} of a functor $\xB\colon\cI\to\cB$ with respect to~$\Phi$ is a functor $\xA\colon\cI\to\cA$ such that~$\Phi\circ\xA$ is naturally equivalent to~$\xB$. In particular, in case~$\cI$ is the one-object, one-morphism category, we identify the functors from~$\cI$ to any category~$\cC$ with the objects of~$\cC$, so a lifting of an object~$B$ of~$\cB$ is an object~$A$ of~$\cA$ such that $\Phi(A)\cong B$.

Our examples below will involve the following categories:

\begin{itemize}
\item The category $\DSLat$ (resp., $\DSLatemb$) of all distributive \jzs s with \jzh s (resp., \jze s).

\item The category $\DLat$ (resp., $\DLatemb$) of all distributive $0$-lattices with $0$-lattice homomorphisms (resp., $0$-lattice embeddings).

\item The category $\Lat$ of all lattices with lattice homomorphisms.
\end{itemize}

Prominent results of lifting functors with respect to the functor $\Conc\colon\Lat\to\DSLat$ are the following:
\begin{itemize}
\item[(1)] E.\,T. Schmidt proved \cite{Schm81} that every distributive $0$-lattice is isomorphic to $\Conc L$ for some lattice~$L$.

\item[(2)] Schmidt's result got extended in 1985 by P. Pudl\'ak \cite{Pudl}, who proved that the inclusion functor $\DLatemb\into\DSLat$ has a lifting $\xA\colon\DLatemb\to\Lat$ with respect to the $\Conc$ functor. So $\Conc\xA(D)\cong D$ naturally in~$D$, for every distributive $0$-lattice~$D$. Furthermore, in Pudl\'ak's construction, $\xA(D)$ is a finite atomistic lattice whenever~$D$ is finite.

\item[(3)] Pudl\'ak's result got further extended by P. R\r u\v zi\v cka \cite{Ruz1}, with a different, ring-theoretical construction that implies that $\xA(D)$ can be taken locally finite, sectionally complemented, and modular, for every distributive $0$-lattice~$D$.

\item[(4)] On the negative side, Pudl\'ak conjectured in 1985 the existence of a lifting of the inclusion functor $\DSLatemb\into\DSLat$ with respect to the~$\Conc$ functor. This conjecture got disproved, before the final negative solution for CLP was obtained,  by J.~T\r uma and F.~Wehrung \cite{Bowtie}.

\end{itemize}

We shall often identify every poset~$K$ with the category whose objects are the elements of~$K$ and where there exists at most one morphism from~$x$ to~$y$, for elements $x,y\in K$, and this occurs exactly in case $x\leq y$.
Denote by $\VPmeas$ the category whose objects are all triples $\seq{P,\mu,S}$, where $P$ is a \mzs, $S$ is a distributive \jzs, and $\mu\colon P\times P\to S$ is a p-measure satisfying the conditions (P1)--(P3) stated in the Abstract, and where the morphisms from $\seq{P,\mu,S}$ to $\seq{Q,\nu,T}$ are the pairs $\seq{f,\boldsymbol{f}}$, where $f\colon P\to Q$ is order-preserving, $\boldsymbol{f}\colon S\to T$ is a \jzh, and $\nu(f(x),f(y))=\boldsymbol{f}(\mu(x,y))$ for all $x,y\in P$. The main result of the present paper (Theorem~\ref{T:PosetRepr}) implies that every functor $\vec{S}\colon K\to\DSLatemb$, where $K$ is (the category associated to) a lattice, has a lifting with respect to the forgetful functor $\Pi\colon \VPmeas\to\DSLat$. (By using the results of~\cite{Ultrabool,RetrLift}, this can be extended to functors $\vec{S}\colon K\to\DSLat$, still for a lattice~$K$, but we shall not present more details about this here.)
Hence, to every distributive \jzs\ $S$, this lifting associates, in a somewhat `natural' fashion, an object of $\VPmeas $ of the form $\seq{P,\mu,S}$.

\subsection{Basic notation and terminology}\label{Su:Basic}

For elements $a$ and $b$ in a poset $P$, we shall use the
abbreviations
 \begin{align*}
 a\prec_P b&\ \Longleftrightarrow\ (a<_Pb
 \text{ and there is no }x
 \text{ such that }a<_Px<_Pb);\\
 a\preceq_P b&\ \Longleftrightarrow\ 
 (\text{either }a\prec_P b\text{ or }a=b);\\
 a\sim_P b&\ \Longleftrightarrow\
 (\text{either }a\leq_Pb\text{ or }b\leq_Pa);\\
 a\inc_P b&\ \Longleftrightarrow\ 
 (a\nleq_Pb\text{ and }b\nleq_Pa).
 \end{align*}
We shall use $\leq$ (instead of $\leq_P$), $\prec$, $\preceq$,
$\sim$, or $\inc$ in case $P$ is understood. We say that~$P$ is \emph{lower finite}, if the principal ideal $\setm{x\in P}{x\leq a}$ is
finite for all $a\in\Lambda$. Observe that in case~$P$ is a \ms, this implies that~$P$ has a least element.

For a category~$\cC$ and a poset~$P$, a \emph{$P$-indexed diagram} in
$\cC$ is a functor $\vec{D}\colon P\to\cC$ (where~$P$ is identified with the associated category). This amounts to a family $\seqm{D_x}{x\in P}$ of objects of~$\cC$, together with
a system of morphisms $\varphi_{x,y}\colon D_x\to D_y$,
for $x\leq y$ in $P$, such that $\varphi_{x,x}=\id_{D_x}$
and $\varphi_{x,z}=\varphi_{y,z}\circ\varphi_{x,y}$ for all
$x\leq y\leq z$ in~$P$, and then we write $\vec{D}=\seqm{D_x,\varphi_{x,y}}{x\leq y\text{ in }P}$.
We shall also denote by
$\vec{D}\res_{\leq p}$ (resp., $\vec{D}\res_{<p}$) the restriction
of $\vec{D}$ to $\setm{x\in P}{x\leq p}$ (resp.,
$\setm{x\in P}{x<p}$), for all $p\in P$.

A \js\ $S$ is \emph{distributive}, if for all $\ba,\bb,\bc\in S$,
if $\bc\leq\ba\vee\bb$, then there are $\bx\leq\ba$ and
$\by\leq\bb$ in $S$ such that $\bc=\bx\vee\by$. Equivalently,
the ideal lattice of $S$ is a distributive lattice, see
\cite[Section~II.5]{GLT2}.

We shall identify every natural number $n$ with the set
$\set{0,1,\dots,n-1}$. We shall denote by $\Pow(X)$ the powerset
of a set~$X$.

\section{Structure of the proof}\label{S:Structure}

\subsection{First obstacle: there is no sequential proof}\label{Su:NoSeq}
Many proofs of positive representation results use transfinite
iterations of `one-step constructions', each of them adding a small number of elements at a time. This is typically the case for J\'onsson's proof of Whitman's Embedding Theorem (cf.~\cite{Jons}). Other examples are the main construction of \cite{FPLat} (that proves,
among other things, that every lattice $L$ such that $\Conc L$ is
a lattice admits a relatively complemented congruence-preserving
extension), or the construction used in~\cite{RTW} to prove that
every distributive
\jzs\ is the range of some `V-distance' of type~$2$, or the
construction used in~\cite{Lamp82} to establish that every
algebraic lattice with compact unit is isomorphic to the congruence
lattice of some groupoid.

However, our result cannot be proved
in such a way. The reason for this is contained
in~\cite{NonExt}, where we construct, \emph{at poset level}, an
example of a p-measure that cannot be extended to a p-measure satisfying~(P2).
This partly explains the complexity of our main construction:
the posets and measures require a somehow `explicit'
construction, which in turn requires quite a large technical
background.

\subsection{General principle of the proof}\label{Su:General}
We need to lift a given \jzs~$S$ with respect to the functor~$\Pi$ introduced in Subsection~\ref{Su:Lifting}. If this is done in case~$S$ has a largest element, the general result follows easily from restricting any p-measure $\mu\colon P\times P\to S\cup\set{1}$ representing~$S\cup\set{1}$ to a suitable lower subset of~$P$ (cf. proof of Corollary~\ref{C:PosetRepr}).

So suppose, from now on, that~$S$ is a \jzus. By Lemma~\ref{L:ErsPudl}, $S$ is the directed union of a family $\vec{D}=\seqm{D_i}{i\in\Lambda}$ of finite distributive \jzu-subsemilattices; furthermore, $\Lambda $ can be taken the collection of all finite subsets of~$S$, in particular \emph{$\Lambda $ is a lower finite lattice}. Our proof will construct a lifting, with respect to the functor~$\Pi$, of~$\vec{D}$ (viewed as a $\Lambda$-indexed diagram); the representation result for~$S$ will follow immediately (cf. proof of Corollary~\ref{C:PosetRepr}).

So now we start with a lower finite \ms\ $\Lambda$ and a $\Lambda$-indexed diagram $\vec{D}=\seqm{D_i,\varphi_{i,j}}{i\leq j\text{ in }\Lambda}$ of finite distributive lattices and \jzuh s. (No stage of the proof will require the totality of these assumptions, nevertheless we shall assume them altogether in the present outline.) As~$\Lambda$ is lower finite, it is well-founded. Accordingly, our lifting will be constructed inductively: from a lifting (with respect to the~$\Pi$ functor) of $\vec{D}\res_{<\ell}$, we shall construct a lifting of $\vec{D}\res_{\leq\ell}$, for any~$\ell\in\Lambda$.

When trying to do this we stumble on a major problem. For such an extension of liftings to be possible, we need a number of somewhat complex assumptions on the lifting of~$\vec{D}\res_{<\ell}$ we are starting with. Fundamental examples, notably in \cite{TuWe1,NonExt}, show that these assumptions cannot be dispensed with.

\subsection{Structure of the lifting at poset level; interval extensions}\label{Su:IntExt}

A first, mild assumption is to suppose that our lifting of~$\vec{D}\res_{<\ell}$ consists, \emph{at poset level}, of inclusion maps. Hence we start with a lifting of~$\vec{D}\res_{<\ell}$ of the form $\seqm{\bQ_i}{i<\ell}$, where $\bQ_i$ is a $D_i$-measured poset for all $i<\ell$ and~$\bQ_j$ is an extension of~$\bQ_i$ with respect to~$\varphi_{i,j}$ for all $i\leq j<\ell$: the latter condition means that~$Q_i$ is a sub-poset of~$Q_j$ and $\bvi{x}{y}_{\bQ_i}=\varphi_{i,j}(\bvi{x}{y}_{\bQ_j})$ for all $x,y\in Q_i$.

Trying to figure out what the poset~$Q_\ell$ should be, an obvious requirement is that it should contain the \emph{set-theoretical union} $P_\ell=\bigcup\famm{Q_i}{i<\ell}$. As~$Q_i$ should be a sub-poset of~$P_\ell$ for each $i<\ell$, we should ensure that the reflexive, transitive binary relation on the set~$P_\ell$ generated by the union of all the partial orderings on the $Q_i$s is a partial ordering (this amounts to verifying that it is antisymmetric).

Easy examples show that this is not the case as a rule.

Hence we need to put conditions on the inclusion maps $Q_i\into Q_j$, for $i\leq j<\ell$. These embeddings will be required to be so-called \emph{interval extensions} (cf. Definition~\ref{D:IntExt}). Furthermore, all the posets~$Q_i$ will be \emph{finite lattices}. As a consequence of the definition of an interval extension, it will turn out that~$Q_i$ is a sublattice of~$Q_j$ for all $i\leq j<\ell$ (cf. Lemma~\ref{L:Abswv}).

\subsection{Normal interval diagrams and covering extensions}\label{Su:CovExt}

In order to endow the set-theoretical union $P_\ell=\bigcup\famm{Q_i}{i<\ell}$ with a suitable poset structure, we need to add, to the condition that~$Q_j$ is an interval extension of~$Q_i$ whenever $i\leq j<\ell$, an assumption of `coherence' between the orderings of the~$Q_i$s. This condition, formulated in Definition~\ref{D:SAsys}, implies, in particular, that $Q_i\cap Q_j=Q_{i\wedge j}$ for all $i,j<\ell$ (remember that~$\Lambda$ is a \ms). We shall say that $\seqm{Q_i}{i<\ell}$ is a \emph{normal interval diagram} of posets. Once this is assumed, the partial ordering on~$P_\ell$ can be easily described from the individual orderings of the~$Q_i$s (cf. Lemma~\ref{L:SAwd}). We shall call~$P_\ell$ the \emph{strong amalgam} of $\seqm{Q_i}{i<\ell}$ (cf. Lemmas~\ref{L:SAwd} and~\ref{L:leqOrd}). It will turn out that~$P_\ell$ is an interval extension of each~$Q_i$ (cf. Lemma~\ref{L:PresRCINT}) and that it is a lattice (cf. Proposition~\ref{P:SA2forP}).

Due to problems pertaining to Example~\ref{Ex:PosetVersDc} and originating in the main counterexample of~\cite{TuWe1}, the p-measures on the~$Q_i$s need not have a common extension to a p-measure on~$P_\ell$. This particular problem is, actually, the hardest technical problem that we need to solve, so we will postpone the required outline until Subsection~\ref{Su:Doubling}.

For the moment, suppose that we have succeeded in finding a p-measure on~$P_\ell$ extending all p-measures on the~$Q_i$s. Even in case the latter p-measures all satisfy Conditions~(P1)--(P3) stated in the Abstract, this may not be the case for $\bvi{{}_-}{{}_-}_{\bP_\ell}$. This is relatively easy to fix, by extending~$P_\ell$ to a larger poset~$Q_\ell$ with a natural extension $\bvi{{}_-}{{}_-}_{\bQ_\ell}$ of $\bvi{{}_-}{{}_-}_{\bP_\ell}$. The general principle underlying the corresponding extension of p-measures is presented in Section~\ref{S:ExtMeasInt}. The poset~$Q_\ell$ will be obtained by inserting the ordinal sum of two suitable finite Boolean lattices in each prime interval of~$P_\ell$ (cf. Proof of Theorem~\ref{T:PosetRepr}). In particular, \emph{$Q_\ell$ is an interval extension of~$P_\ell$}.

At this point, we stumble on the slightly annoying point that an interval extension of an interval extension may not be an interval extension (cf. Example~\ref{Ex:SpNotComp}). So nothing would guarantee \emph{a priori} that~$Q_\ell$ is an interval extension of each~$Q_i$ for $i<\ell$, and the induction process would break down. Fortunately, $Q_\ell$ is what we call in Definition~\ref{D:covering} a \emph{covering extension} of~$P_\ell$, which, by Proposition~\ref{P:CompIntCov}, will be sufficient to ensure that $Q_\ell$ is, indeed, an interval extension of each~$Q_i$. (It would be too much asking that $Q_\ell$ be a covering extension of each~$Q_i$.) Therefore, the induction ball keeps rolling---at least at poset level.

\subsection{The elements $x_{\bl}$ and $x^{\bl}$}\label{Su:xblbl}
A fundamental tool in the inductive evaluation of the Boolean values $\bvi{x}{y}_{\bP_\ell}$ is introduced in Section~\ref{S:x**}. For each $x\in P_\ell\setminus Q_0$, the properties of normal interval diagrams imply the existence of a least element of~$P_\ell$, denoted by $x^\bl$, such that $x<x^\bl$ and~$x^\bl$ lies in a block~$Q_i$ of smaller index than the one containing~$x$ (the latter condition is formulated $\nu(x^\bl)<\nu(x)$ in Lemma~\ref{L:Existx**}---the `complexity' $\nu(x)$ of~$x$ is the least~$i<\ell$ such that $x\in Q_i$). The element~$x_\bl$ is defined dually, so~$x_\bl$ is the largest element smaller than~$x$ such that $\nu(x_\bl)<\nu(x)$. The inductive definition of $\bvi{x}{y}_{\bP_\ell}$ will make a heavy use of the elements~$x^\bl$ and~$x_\bl$.

\subsection{Finding the p-measure on~$P_\ell$: doubling extensions}\label{Su:Doubling}
As mentioned earlier, the hardest technical problem of the whole paper is to ensure that the p-measures on the~$\bQ_i$, for $i<\ell$, have a common extension to some p-measure on~$\bP_\ell$. Recall that Example~\ref{Ex:PosetVersDc} shows that this cannot be done without additional assumptions.

The idea that we implement here is to relate the Boolean value~$\bvi{x}{y}$, for $x,y\in P_\ell$, to Boolean values involving less complicated elements, such as $\bvi{x^\bl}{y}$, $\bvi{x}{y_\bl}$, and so on. The first doubling condition~(DB1) (cf. Section~\ref{S:Doubling}) says that every $x\in P_\ell$ is `closest' (with respect to Boolean values) either to~$x^\bl$ or to~$x_\bl$. The second doubling condition~(DB2) says that if~$x$ is closest to~$x^\bl$ and in `non-degenerate' cases, $\bvi{x}{y}=\bvi{x^\bl}{y}$ whenever both~$x$ and~$y$ (and thus also~$x^\bl$) belong to the same~$Q_i$; and, symmetrically, if~$x$ is closest to~$x_\bl$ and in non-degenerate cases, $\bvi{y}{x}=\bvi{y}{x_\bl}$.

This is the basis for our definition of~$\bvi{x}{y}$, for $x,y\in P_\ell$. In case $x,y\in Q_i$ for some $i<\ell$, we put $\bvo{x}{y}=\varphi_{i,\ell}(\bvi{x}{y}_{\bQ_i})$ (cf.~\eqref{Eq:Defbvoxy}). This is the easiest case where we can evaluate $\bvi{x}{y}$, for $x,y\in P_\ell$---namely by setting $\bvi{x}{y}=\bvo{x}{y}$.

In the general case, the most natural guess is then to define $\bvi{x}{y}$ as the meet, in~$D_\ell$, of all joins of the form
$\bigvee_{i<n}\bvo{z_i}{z_{i+1}}$, where~$n$ is a positive integer, $z_0$, $z_1$, \dots, $z_n$ are elements of~$P_\ell$, $z_0=x$, and $z_n=y$.

Most of the technical difficulties of the paper are contained in the proof that this guess is sound.

This proof makes a heavy use of the doubling conditions. Furthermore, as an unexpected side issue, it implies that it is sufficient to consider the case where $n=3$ (cf. Corollary~\ref{C:Simplebvd}). The propagation of~(DB1) to the level~$\ell$ is taken care of by Lemma~\ref{L:PdoublesQi}. The propagation of~(DB2) is taken care of by Lemma~\ref{L:PropagADD3}.

This takes care of the construction of~$\bvi{{}_-}{{}_-}_{\bP_\ell}$. As mentioned earlier, the technical background for the extension of that p-measure to a suitable p-measure on~$Q_\ell$ is contained in Section~\ref{S:ExtMeasInt}. The construction of the poset~$Q_\ell$ itself is quite natural, and it is contained in the proof of Theorem~\ref{T:PosetRepr}.

Corollary~\ref{C:PosetRepr} trivially implies the result stated in the Abstract.

Further comments, in particular about the possible uses of our main result to tackle current open problems, are presented in Section~\ref {S:Dist}.

\section{Relatively complete and interval extensions of posets}
\label{S:IntExt}

\begin{definition}\label{D:IntExt}
A poset $Q$ is a \emph{relatively complete extension} of a poset
$P$, in notation $P\lerc Q$, if for all $x\in Q$, there exists a
largest element of $P$ below~$x$ (denoted by $x_P$) and a least
element of $P$ above~$x$ (denoted by $x^P$). Then we define
binary relations $\ll_P$ and $\equiv_P$ on $Q$ by
 \begin{align*}
 x\ll_Py&\ \Longleftrightarrow\ x^P\leq y_P,\\
 x\equiv_Py&\ \Longleftrightarrow\ (x_P=y_P\text{ and }x^P=y^P),
 \end{align*}
for all $x,y\in Q$. We say that $Q$ is an \emph{interval
extension} of $P$, in notation $P\leint Q$, if $P\lerc Q$ and
for all $x,y\in Q$, $x\leq y$ implies that either $x\ll_Py$ or
$x\equiv_Py$.
\end{definition}

The proofs of the following two lemmas are easy exercises.

\begin{lemma}\label{L:RCcompose}
Let $P$, $Q$, and $R$ be posets. If $P\lerc Q$ and $Q\lerc R$,
then $P\lerc R$.
\end{lemma}

\begin{lemma}\label{L:Abswv}
Let $Q$ be a relatively complete extension of a poset $P$ and let
$X\subseteq P$. Then $\bigvee_PX$ exists if{f} $\bigvee_QX$
exists, and then the two values are equal. The dual statement
also holds.
\end{lemma}

In particular, if $P\lerc Q$ and $Q$ is a lattice, then $P$ is a
sublattice of $Q$. Hence, from now on, when dealing with
relatively complete extensions, we shall often omit to mention in
which subset the meets and joins are evaluated.

\begin{lemma}\label{L:meetincomp}
Let $Q$ be an interval extension of a poset $P$, let $x\in P$
and $y,z\in Q$. If $x,y\leq z$ and $x\nleq y$, then $y^P\leq z$;
and dually.
\end{lemma}

\begin{proof}
If $y\equiv_Pz$, then, as $x\in P$ and $x\leq z$, we get
$x\leq z_P=y_P\leq y$, \contr; hence
$y\not\equiv_Pz$. As $y\leq z$ and $P\leint Q$, we get
$y\ll_Pz$, and thus $y^P\leq z$.
\end{proof}

\begin{lemma}\label{L:IntLatt}
Let $Q$ be an interval extension of a poset $P$. Then $Q$ is a
lattice if{f}~$P$ is a lattice and the interval $[x_P,x^P]$ is a
lattice for each $x\in Q$. Furthermore, if~$Q$ is a lattice, then
for all incomparable $x,y\in Q$,
 \begin{align}
 x\vee y&=\begin{cases}
 x^P\vee y^P,&\text{if }x\not\equiv_Py,\\
 x\vee_{[u,v]}y,&\text{if }x_P=y_P=u\text{ and }x^P=y^P=v,
 \end{cases}\label{Eq:JoinInt}\\
 x\wedge y&=\begin{cases}
 x^P\wedge y^P,&\text{if }x\not\equiv_Py,\\
 x\wedge_{[u,v]}y,&\text{if }x_P=y_P=u\text{ and }x^P=y^P=v.
 \end{cases}\label{Eq:MeetInt}
 \end{align}
\end{lemma}

\begin{proof}
We prove the nontrivial direction. So suppose that $P$ is a
lattice and the interval $[x_P,x^P]$ is a lattice for each
$x\in Q$. For incomparable $x,y\in P$, we prove that the join
$x\vee y$ is defined in $Q$ and given by \eqref{Eq:JoinInt}. The
proof for the meet is dual. So let $z\in Q$ such that
$x,y\leq z$. If $x\ll_Pz$, then $x^P\leq z$, hence, using
Lemma~\ref{L:meetincomp}~(with $x^P$ instead of~$x$), we obtain
$y^P\leq z$, and hence, using Lemma~\ref{L:Abswv},
$x^P\vee y^P\leq z$. The conclusion is similar for $y\ll_Pz$. As
$P\leint Q$, the remaining case is where $x\equiv_Pz\equiv_Py$.
Putting $u=x_P=y_P$ and $v=x^P=y^P$, the interval $[u,v]$ is, by
assumption, a lattice, so $x\vee_{[u,v]}y\leq z$, and hence
$x\vee y=x\vee_{[u,v]}y\leq z$.
\end{proof}

\begin{lemma}\label{L:compxxP}
Let $Q$ be an interval extension of a poset $P$. Then $x\sim y$
implies that $x_P\sim y$ and $x^P\sim y$, for all $x,y\in Q$.
\end{lemma}

\begin{proof}
We prove the result for $x_P$. If $x\leq y$, then $x_P\leq y$
and we are done. Suppose that $y\leq x$. If $y\ll_Px$, then
$y\leq x_P$. Suppose that $y\not\ll_Px$. As
$y\leq x$ and $P\leint Q$, we get $x\equiv_Py$, and thus
$x_P=y_P\leq y$.
\end{proof}

\begin{definition}\label{D:SISch}
A \emph{standard interval scheme} is a family of the form
\linebreak
$\seq{P,\seqm{Q_{a,b}}{\seq{a,b}\in\mathcal{I}}}$, where the
following conditions are satisfied:
\begin{enumerate}
\item $P$ is a poset, $\mathcal{I}$ is a subset of
$\setm{\seq{a,b}\in P\times P}{a<b}$, and $Q_{a,b}$ is a
(possibly empty) poset, for all $\seq{a,b}\in\mathcal{I}$.

\item $Q_{a,b}\cap P=\es$ for all
$\seq{a,b}\in\mathcal{I}$.

\item $Q_{a,b}\cap Q_{c,d}=\es$ for all distinct
$\seq{a,b},\seq{c,d}\in\mathcal{I}$.
\end{enumerate}
We say that the standard interval scheme above is \emph{based on
$P$}.
\end{definition}

The proofs of the following two lemmas are straightforward
exercises.

\begin{lemma}\label{L:SumSISch}
Let $\seq{P,\seqm{Q_{a,b}}{\seq{a,b}\in\mathcal{I}}}$ be a
standard interval scheme. Put
$Q=P\cup\bigcup\famm{Q_{a,b}}{\seq{a,b}\in\mathcal{I}}$.
Furthermore, for all $x\in Q$, put $x_P=x^P=x$ if $x\in P$,
while $x_P=a$ and $x^P=b$ if $x\in Q_{a,b}$, for
$\seq{a,b}\in\mathcal{I}$. Let $x\leq y$ hold, if either
$x^P\leq y_P$ or there exists $\seq{a,b}\in\mathcal{I}$ such
that $x,y\in Q_{a,b}$ and $x\leq_{Q_{a,b}}y$, for all $x,y\in Q$.
Then $\leq$ is a partial ordering on $Q$ and $Q$ is an
interval extension of $P$.
\end{lemma}

In the context of Lemma~\ref{L:SumSISch}, we shall use the
notation
 \begin{equation}\label{Eq:SumNotSISch}
 Q=P+\sum\famm{Q_{a,b}}{\seq{a,b}\in\mathcal{I}}.
 \end{equation}
Conversely, the following lemma shows that any interval extension
can be obtained by the
$P+\sum\famm{Q_{a,b}}{\seq{a,b}\in\mathcal{I}}$ construction.
This construction is a special case of a construction presented
in \cite{GrKe}.

\begin{lemma}\label{L:FromInt2SIS}
Let $Q$ be an interval extension of a poset $P$. Put\linebreak
$\mathcal{I}=\setm{\seq{a,b}\in P\times P}{a<b}$, and
$Q_{a,b}=\setm{x\in Q}{x_P=a\text{ and }x^P=b}$,
for all $\seq{a,b}\in P$. Then
$\seq{P,\seqm{Q_{a,b}}{\seq{a,b}\in\mathcal{I}}}$ is a standard
interval scheme, and
$Q=P+\sum\famm{Q_{a,b}}{\seq{a,b}\in\mathcal{I}}$.
\end{lemma}

It follows from Lemmas~\ref{L:SumSISch}
and~\ref{L:FromInt2SIS} that any standard interval scheme based
on~$P$ defines an interval extension of $P$, and every interval
extension of $P$ is defined \emph{via} some standard interval
scheme on $P$.

\section{Covering extensions of posets}\label{S:Cov}

\begin{definition}\label{D:covering}
We say that a poset $Q$ is a \emph{covering extension} of a
poset $P$, in notation $P\lecov Q$, if $Q$ is an interval
extension of $P$ (cf. Definition~\ref{D:IntExt}) and the relation
$x_P\preceq_Px^P$ holds for all $x\in Q$.
\end{definition}

\begin{lemma}\label{L:Int+Cov2Int}
Let $P$, $Q$, and $R$ be posets such that $P\leint Q$,
$Q\leint R$, and there are no $x\in P$ and $y\in R$ such
that $y_Q<x<y^Q$. Then $P\leint R$.
\end{lemma}

\begin{proof}
First, $P\lerc R$ (cf. Lemma~\ref{L:RCcompose}). Now let
$x\leq y$ in $R$, and assume, towards a contradiction, that
$x\not\ll_Py$ and $x\not\equiv_Py$.

If $x\equiv_Qy$, then $x\equiv_Py$, \contr. As $x\leq y$ and
$Q\leint R$, we get $x\ll_Qy$, that is, $x^Q\leq y_Q$. If
$x^Q\ll_Py_Q$, then $x\ll_Py$, \contr. As $x^Q\leq y_Q$ and
$P\leint Q$, we get $x^Q\equiv_Py_Q$.

As $y_Q\leq y^Q$ and $P\leint Q$, either $y_Q\equiv_Py^Q$ or
$y_Q\ll_Py^Q$. In the first case, we get, using the relation
$x^Q\equiv_Py_Q$, the equalities
$x^P=(x^Q)^P=(y_Q)^P=(y^Q)^P=y^P$. In the second case, we get,
using  again the relation $x^Q\equiv_Py_Q$, the
inequalities $x^P=(x^Q)^P=(y_Q)^P\leq y^Q$. But
$y_Q\leq(y_Q)^P=x^P$, and so $y_Q\leq x^P\leq y^Q$. Hence, by
assumption, either $x^P=y_Q$ or $x^P=y^Q$. If $x^P=y_Q$, then
$x^P=y_P$, so $x\ll_Py$, \contr; hence only the subcase
where $x^P=y^Q$ remains, so $y^Q\in P$, and so $x^P=y^Q=y^P$.

So we have proved that in either case, the equality $x^P=y^P$
holds. Dually, the equality $x_P=y_P$ holds, and so
$x\equiv_Py$, \contr.
\end{proof}

\begin{proposition}\label{P:CompIntCov}
For arbitrary posets $P$, $Q$, and $R$, the following statements
hold:
\begin{enumerate}
\item If $P\leint Q$ and $Q\lecov R$, then $P\leint R$.

\item If $P\lecov Q$ and $Q\lecov R$, then $P\lecov R$.
\end{enumerate}
\end{proposition}

\begin{proof}
(i) follows immediately from Lemma~\ref{L:Int+Cov2Int}. Now we
prove (ii). So assume that $P\lecov Q$ and $Q\lecov R$, let
$x\in R$, we prove that $x_P\preceq_Px^P$. If
$\set{x_Q,x^Q}\subseteq P$, then $x_P=x_Q$ and $x^P=x^Q$, but
$Q\lecov R$, thus $x_Q\preceq_Qx^Q$, and thus, \emph{a fortiori},
$x_P\preceq_Px^P$. So suppose, from now on, that
$\set{x_Q,x^Q}\not\subseteq P$, say $x^Q\notin P$.

If $(x^Q)_P\nleq x$, then, as $x,(x^Q)_P\leq x^Q$, as $P\leint R$
(proved in (i)), and by Lemma~\ref{L:meetincomp}, we get
$x^P\leq x^Q$, so $x^Q\in P$, \contr. Hence $(x^Q)_P\leq x$, but
$(x^Q)_P$ lies above $x_P$ and belongs to $P$, and so
$(x^Q)_P=x_P$. As $P\lecov Q$, we get
$x_P=(x^Q)_P\preceq_P(x^Q)^P=x^P$.
\end{proof}

\begin{example}\label{Ex:SpNotComp}
The following example shows that interval extensions do not
compose. We let $K$, $L$, and $M$ the lattices diagrammed on
Figure~\ref{Fig:SpNotComp}. Then $K\leint L$ and $L\leint M$,
however $K\not\leint M$, as $x\leq y$ and
$x\not\equiv_Ky$ while $x^K\nleq y_K$. In this example,
$K\lecov L$ and $L\not\lecov M$.

\begin{figure}[htb]
\includegraphics{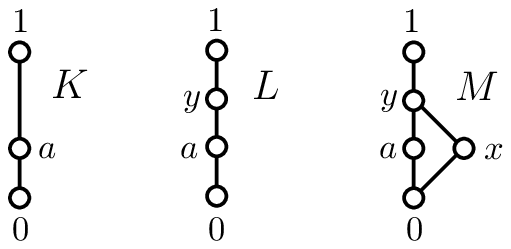}
\caption{Interval extensions do not compose.}
\label{Fig:SpNotComp}
\end{figure}
\end{example}

\section{Strong amalgams of normal diagrams of posets}
\label{S:SASIS}

In this section we shall deal with families of posets indexed by
\ms s.

\begin{definition}\label{D:SAsys}
A \emph{normal diagram of posets} consists
of a family $\vec{Q}=\seqm{Q_i}{i\in\Lambda}$ of posets, indexed
by a \ms\ $\Lambda$, such that the following conditions hold (we
denote by $\leq_i$ the partial ordering of $Q_i$, for all
$i\in\Lambda$):
\begin{itemize}
\item[(1)] $Q_i$ is a sub-poset of $Q_j$ for all $i\leq j$ in
$\Lambda$.

\item[(2)] $Q_i\cap Q_j=Q_{i\wedge j}$
(set-theoretically!) for all $i,j\in\Lambda$.

\item[(3)] For all $i,j,k\in\Lambda$ such that $i,j\leq k$ and
all $\seq{x,y}\in Q_i\times Q_j$, if $x\leq_k y$, then there
exists $z\in Q_{i\wedge j}$ such that $x\leq_iz$ and $z\leq_jy$.
\end{itemize}
Furthermore, we say that $\vec{Q}$ is a
\emph{normal interval diagram of posets}, if $Q_j$ is an
interval extension of $Q_i$ for all~$i\leq j$ in~$\Lambda$.
\end{definition}

Let $\vec{Q}=\seqm{Q_i}{i\in\Lambda}$ be a normal diagram of
posets and set $P=\bigcup\famm{Q_i}{i\in\Lambda}$. For
$x,y\in P$, let $x\leq y$ hold, if there are
$i,j\in\Lambda$ and $z\in Q_{i\wedge j}$ such that
$x\in\nobreak Q_i$, $y\in\nobreak Q_j$, and $x\leq_iz\leq_jy$. As
the following lemma shows, this definition is independent of the
chosen pair $\seq{i,j}$ such that $x\in Q_i$ and $y\in Q_j$.

\begin{lemma}\label{L:SAwd}
For all $x,y\in P$ and all $i,j\in\Lambda$ such that $x\in Q_i$
and $y\in Q_j$, $x\leq y$ if{f} there exists $z\in Q_{i\wedge j}$
such that $x\leq_iz\leq_jy$.
\end{lemma}

\begin{proof}
The given condition implies, by definition, that $x\leq y$.
Conversely, suppose that $x\leq y$, and fix $i',j'\in\Lambda$
and $z'\in Q_{i'\wedge j'}$ such that $x\in Q_{i'}$,
$y\in Q_{j'}$, and $x\leq_{i'}z'\leq_{j'}y$. As
$x\in Q_{i\wedge i'}$, $z'\in Q_{i'\wedge j'}$, and
$x\leq_{i'}z'$, there exists
$t\in Q_{i\wedge i'\wedge j'}$ such that
$x\leq_{i\wedge i'}t\leq_{i'\wedge j'}z'$. As
$t\in Q_{i\wedge i'\wedge j'}$, $y\in Q_{j\wedge j'}$, and
$t\leq_{j'}y$, there exists
$z\in Q_{i\wedge i'\wedge j\wedge j'}$ such that
$t\leq_{i\wedge i'\wedge j'}z\leq_{j\wedge j'}y$. In particular,
$z\in Q_{i\wedge j}$ and $x\leq_iz\leq_jy$.
\end{proof}

\begin{lemma}\label{L:leqOrd}
The binary relation $\leq$ defined above is a partial ordering
of $P$. Furthermore, $Q_i$ is a sub-poset of $P$ for all
$i\in\Lambda$.
\end{lemma}

\begin{proof}
Reflexivity is obvious. Now let $x,y\in P$ such
that $x\leq y$ and $y\leq x$. Fix $i,j\in\Lambda$ such that
$x\in Q_i$ and $y\in Q_j$. By
Lemma~\ref{L:SAwd}, there are $u,v\in Q_{i\wedge j}$ such that
$x\leq_iu\leq_jy\leq_jv\leq_ix$. Hence
$u\leq_{i\wedge j}v\leq_{i\wedge j}u$, and so $u=v$, and
therefore $x=u=y$.

Now let $x\leq y\leq z$ in $P$, and fix $i,j,k\in\Lambda$
such that $x\in Q_i$, $y\in Q_j$, and $z\in Q_k$.
By Lemma~\ref{L:SAwd}, there are $u\in Q_{i\wedge j}$ and
$v\in Q_{j\wedge k}$ such that $x\leq_iu\leq_jy$ and
$y\leq_jv\leq_kz$. As $u\leq_jv$, there exists
$w\in Q_{i\wedge j\wedge k}$ such that
$u\leq_{i\wedge j}w\leq_{j\wedge k}v$. Hence
$x\leq_iw\leq_kz$, and so $x\leq z$. Therefore, $\leq$ is a
partial ordering on $P$.

Finally, let $i\in\Lambda$ and let $x,y\in Q_i$. If
$x\leq_iy$, then $x\leq y$ trivially. Conversely, if $x\leq y$,
then, by Lemma~\ref{L:SAwd}, there exists $z\in Q_i$ such that
$x\leq_iz\leq_iy$, whence $x\leq_iy$. Therefore, $x\leq y$
if{f} $x\leq_iy$.
\end{proof}

Hence, from now on, we shall most of the time drop the index $i$ in $\leq_i$,
for $i\in\Lambda$. We shall call the poset~$P$ the \emph{strong
amalgam} of $\seqm{Q_i}{i\in\Lambda}$.

\begin{lemma}\label{L:PresRCINT}
Let $\vec{Q}=\seqm{Q_i}{i\in\Lambda}$ be a normal diagram of
posets, with strong amalgam $P=\bigcup\famm{Q_i}{i\in\Lambda}$.
Then the following statements hold:
\begin{enumerate}
\item If $i\leq j$ implies that $Q_i\lerc Q_j$ for all
$i,j\in\Lambda$, then $Q_i\lerc P$ for all $i\in\Lambda$.

\item If $i\leq j$ implies that $Q_i\leint Q_j$ for all
$i,j\in\Lambda$, then $Q_i\leint P$ for all $i\in\Lambda$.
\end{enumerate}
\end{lemma}

\begin{proof}
(i). Let $x\in P$, say $x\in Q_j$, for $j\in\Lambda$,
and let $i\in\Lambda$. By Lemma~\ref{L:SAwd}, every
element of~$Q_i$ below~$x$ lies below $x_{Q_{i\wedge j}}$; hence
$x_{Q_i}$ exists, and it is equal to $x_{Q_{i\wedge j}}$.
Dually, $x^{Q_i}$ exists, and it is equal to $x^{Q_{i\wedge j}}$.
In particular, $Q_i\lerc P$.

\begin{quote}\em
To ease notation, we shall from now on use the abbreviations
$\lc{x}{i}=x_{Q_i}$ and $\uc{x}{i}=x^{Q_i}$, for all $x\in P$
and all $i\in\Lambda$. Similarly, we shall abbreviate
$x\equiv_{Q_i}y$ by $x\equiv_iy$ and $x\ll_{Q_i}y$ by $x\ll_iy$,
for all $x,y\in P$ and all $i\in\Lambda$.
\end{quote}

(ii). First, it follows from (i) above that $Q_i\lerc P$.
Now let $x,y\in P$ such that $x\leq y$, we prove that either
$x\equiv_iy$ or $x\ll_iy$. Fix $j,k\in\Lambda$ such that
$x\in Q_j$ and $y\in Q_k$. Suppose first that $j=k$. As
$Q_{i\wedge j}\leint Q_j$, either $x\equiv_{i\wedge j}y$ or
$x\ll_{i\wedge j}y$. As $\lc{t}{i}=\lc{t}{i\wedge j}$ and
$\uc{t}{i}=\uc{t}{i\wedge j}$ for all $t\in\set{x,y}$
(see proof of (i) above), this
amounts to saying that either $x\equiv_iy$ or $x\ll_iy$, so we
are done.

In the general case, there exists, by Lemma~\ref{L:SAwd},
$z\in Q_{j\wedge k}$ such that $x\leq_jz\leq_ky$. Applying the
paragraph above to the pairs $\seq{x,z}$ and $\seq{z,y}$, we
obtain that either $x\equiv_iz$ or $x\ll_iz$, and either
$z\equiv_iy$ or $z\ll_iy$. If $x\equiv_iz$ and $z\equiv_iy$,
then $x\equiv_iy$. In all other three cases, $x\ll_iy$.
\end{proof}

\begin{proposition}\label{P:SA2forP}
Let $\vec{Q}=\seqm{Q_i}{i\in\Lambda}$ be a normal interval
diagram of \emph{lattices}. Then the following statements hold:
\begin{enumerate}
\item The strong amalgam $P=\bigcup\famm{Q_i}{i\in\Lambda}$ is a
lattice.

\item $Q_i$ is a sublattice of $P$ for all $i\in\Lambda$.

\item For all $i,j\in\Lambda$ and all incomparable $a\in Q_i$
and $b\in Q_j$, both $a\vee b$ and $a\wedge b$ belong to
$Q_{i\wedge j}$.
\end{enumerate}
\end{proposition}

\begin{proof}
We denote by $\vee_k$ (resp., $\wedge_k$) the join (resp., meet)
operation in $Q_k$, for all $k\in\Lambda$.
We first establish a claim.

\begin{sclaim}
Let $i,j,k\in\Lambda$ with $i,j\leq k$ and let
$\seq{x,y}\in Q_i\times Q_j$. If $x\inc y$, then both
$x\vee_ky$ and $x\wedge_ky$ belong to $Q_{i\wedge j}$.
\end{sclaim}

\begin{scproof}
If $x\equiv_iy$, then, as $x\in Q_i$, we obtain that $x=y$,
which contradicts the assumption that $x\inc y$; hence
$x\not\equiv_iy$. As $Q_i\leint Q_k$, it follows from
Lemma~\ref{L:IntLatt} that
$x\vee_ky=\uc{x}{i}\vee_k\uc{y}{i}$, and thus, as $Q_i$ is a
sublattice of $Q_k$, $x\vee_ky\in Q_i$. Similarly,
$x\vee_ky\in Q_j$, and hence $x\vee_ky\in Q_{i\wedge j}$. The
proof for the meet is dual.
\end{scproof}

Now we establish (iii). We give the proof for the meet; the
proof for the join is dual. Suppose that $a\inc b$, let
$i,j\in\Lambda$ such that $a\in Q_i$ and
$b\in Q_j$, and put
$c=\lc{a}{i\wedge j}\wedge_{i\wedge j}\lc{b}{i\wedge j}$. Of
course, $c\leq a,b$. Now let $x\in P$ such that $x\leq a,b$, we
prove that $x\leq c$. Pick $k\in\Lambda$ such that $x\in Q_k$
and set $m=i\wedge j\wedge k$.
By Lemma~\ref{L:SAwd}, there are
$a'\in Q_{i\wedge k}$ and $b'\in Q_{j\wedge k}$ such that
$x\leq a'\leq a$ and $x\leq b'\leq b$. Suppose first that
$a'\inc b'$. It follows from the Claim above that $a'\wedge_kb'$
belongs to~$Q_m$, thus to $Q_{i\wedge j}$. As
$x\leq a'\wedge_kb'\leq a,b$, we obtain that
$x\leq a'\wedge_kb'\leq c$.

Suppose now that $a'\sim b'$, say $a'\leq b'$.
By Lemma~\ref{L:SAwd}, there
exists $a''\in Q_m$ such that $a'\leq a''\leq b'$. If
$a''\leq a$, then (as $a''\leq b$ and $a''\in Q_{i\wedge j}$)
$a''\leq c$, and so $x\leq c$. As $a\nleq a''$ (for $a\nleq b$),
the only possibility left is $a\inc a''$. By the Claim above,
$a\wedge_ia''$ belongs to $Q_m$, but this element lies below
both $a$ and $b$, thus, again, below~$c$. As
$x\leq a'\leq a\wedge_ia''$, we thus obtain that $x\leq c$.

Hence $c$ is the meet of $\set{a,b}$ in $P$, and so $P$ is a \ms.
Dually, $P$ is a \js; this establishes (i).
In case $a,b\in Q_i$, we take $i=j$, and thus $c=a\wedge_ib$, and
so we obtain that $Q_i$ is a meet-subsemilattice of $P$.
Dually, $Q_i$ is a join-subsemilattice of $P$; this
establishes~(ii).
\end{proof}

\section{Lower finite normal interval diagrams; the elements
$x_{\bl}$ and $x^{\bl}$}\label{S:x**}

For a normal diagram $\vec{Q}=\seqm{Q_i}{i\in\Lambda}$ of
posets, it follows from Definition~\ref{D:SAsys}(2)
that for every element
$x$ of $P=\bigcup\famm{Q_i}{i\in\Lambda}$, the set
$\setm{i\in\Lambda}{x\in Q_i}$ is closed under finite meets.
In particular, in case $\Lambda$ is lower finite (cf. Subsection~\ref{Su:Basic}; we shall say that the diagram~\emph{$\vec{Q}$ is lower finite}), there exists
a least $i\in\Lambda$ such that $x\in Q_i$. We shall denote this
element by~$\nu(x)$, and we shall call the map
$\nu\colon P\to\Lambda$ the \emph{valuation} associated with the
normal diagram~$\vec{Q}$.

In this section, we shall fix a lower finite normal
interval diagram $\vec{Q}=\seqm{Q_i}{i\in\Lambda}$ of lattices,
with strong amalgam $P=\bigcup\famm{Q_i}{i\in\Lambda}$. It
follows from Proposition~\ref{P:SA2forP} that $P$ is a lattice.

\begin{lemma}\label{L:Existx**}
For all $x\in P\setminus Q_0$,
there exist a largest $x_{\bl}<x$ such that
$\nu(x_{\bl})<\nu(x)$ and a least $x^{\bl}>x$ such that
$\nu(x^{\bl})<\nu(x)$. Furthermore, the following hold:
\begin{enumerate}
\item $\nu(x_{\bl})$ and $\nu(x^{\bl})$ are comparable.

\item Putting
$i=\max\set{\nu(x_{\bl}),\nu(x^{\bl})}$, both
equalities $x_{\bl}=\lc{x}{i}$ and $x^{\bl}=\uc{x}{i}$ hold.

\item For all $y\in P$ such that $\nu(x)\nleq\nu(y)$, $x\leq y$
implies that $x^{\bl}\leq y$, and $y\leq x$ implies that
$y\leq x_{\bl}$.
\end{enumerate}
\end{lemma}

\begin{proof}
Put $\Lambda'=\setm{i\in\Lambda}{i<\nu(x)}$ and
$X=\setm{\lc{x}{i}}{i\in\Lambda'}$. As $x\notin Q_0$, the
set~$X$ is nonempty. As $X$ is finite (because $\Lambda'$ is
finite), it has a join in $P$, say $x_{\bl}$. It follows easily
from Proposition~\ref{P:SA2forP} that
$\nu(x_{\bl})<\nu(x)$, whence $x_{\bl}$ is the largest element
of~$X$. The proof of the existence of $x^{\bl}$ is similar.

As $x_{\bl}\leq x^{\bl}$, there exists
$y\in Q_{\nu(x_{\bl})\wedge\nu(x^{\bl})}$ such that
$x_{\bl}\leq y\leq x^{\bl}$. If $x\leq y$, then, as
$\nu(y)<\nu(x)$, we get
$y=x^{\bl}$, and thus $\nu(x^{\bl})\leq\nu(x_{\bl})$.
Similarly, if $y\leq x$, then $y=x_{\bl}$, and thus
$\nu(x_{\bl})\leq\nu(x^{\bl})$. Now suppose that
$x\inc y$. It follows from Proposition~\ref{P:SA2forP} that
$\nu(x\wedge y),\nu(x\vee y)\leq\nu(y)<\nu(x)$, thus, as
$x\vee y\leq x^{\bl}$ and $x_{\bl}\leq x\wedge y$, we get
$x\vee y=x^{\bl}$ and $x\wedge y=x_{\bl}$. By using the first
equality, we get
$\nu(x^{\bl})=\nu(x\vee y)\leq\nu(y)\leq\nu(x_{\bl})$, while by
using the second one, we get $\nu(x_{\bl})\leq\nu(x^{\bl})$, and
hence $\nu(x_{\bl})=\nu(x^{\bl})$. This takes care of (i).

Now we deal with (ii). {}From
$\nu(\lc{x}{i}),\nu(\uc{x}{i})\leq i<\nu(x)$ it follows that
$\lc{x}{i}\leq x_{\bl}$ and $x^{\bl}\leq\uc{x}{i}$. As both
$x_{\bl}$ and $x^{\bl}$ belong to $Q_i$, we get
$\lc{x}{i}=x_{\bl}$ and $\uc{x}{i}=x^{\bl}$.

Let $y\in P$ with $\nu(x)\nleq\nu(y)$. If $x\leq y$, then
there exists
$z\in Q_{\nu(x)\wedge\nu(y)}$ such that $x\leq z\leq y$, but
$\nu(z)<\nu(x)$, thus $x^{\bl}\leq z$, and so $x^{\bl}\leq y$.
The proof for $x_{\bl}$ is dual. This takes care of~(iii).
\end{proof}

\section{Extending a p-measure to an interval extension}
\label{S:ExtMeasInt}

Recall that we introduced p-measures and p-measured posets in Definition~\ref{D:PosetMeas}.

We shall define the \emph{distance function} on a p-measured
poset $\bP$ by
 \[
 \bve{x}{y}_{\bP}=\bvP{\max\set{x,y}}{\min\set{x,y}},
 \qquad\text{for all \emph{comparable} }x,y\in P.
 \]
Obviously, the distance function on $P$ satisfies the
\emph{triangular inequality}\linebreak
$\bve{x}{z}_{\bP}\leq\bve{x}{y}_{\bP}\vee\bve{y}{z}_{\bP}$, for
all pairwise comparable $x,y,z\in P$. Furthermore, the equality
holds for $x\geq y\geq z$.

For \jzs s $S$ and $T$ and a \jzh\ $\varphi\colon S\to T$,
a $S$-valued p-measured poset $\bP$, and a $T$-valued p-measured
poset $\bQ$, we shall say that~$\bQ$ \emph{extends}
$\bP$ with respect to $\varphi$, if $P$ is a sub-poset of $Q$ and
 \[
 \bvQ{x}{y}=\varphi\left(\bvP{x}{y}\right),\qquad
 \text{for all }x,y\in P.
 \]
We shall then say that the inclusion map from $P$ into $Q$,
together with $\varphi$, form a \emph{morphism} from $\bP$ to
$\bQ$, and define diagrams of p-measured posets
accordingly. (Obviously, we could have defined morphisms more
generally by involving an order-embedding from $P$ into $Q$, but
the present definition is sufficient, and more convenient, for
our purposes.)

Until Lemma~\ref{L:TrIneqbvQ}, we fix a distributive lattice $D$
with zero and a $D$-valued p-measured poset $\bP$. We
are given an interval extension $Q$ of $P$ in which each
interval of $Q$ of the form $[x_P,x^P]$, for $x\in Q$, is
endowed with a p-measure $\bvi{{}_{-}}{{}_{-}}_{[x_P,x^P]}$.
We assume \emph{compatibility} between those p-measures, in the
sense that $\bve{x^P}{x_P}_{\bP}=\bve{x^P}{x_P}_{[x_P,x^P]}$, for
all $x\in Q$. We define a map $\bvi{{}_{-}}{{}_{-}}_{\bQ}$ from
$Q\times Q$ to~$D$, by setting
$\bvQ{x}{y}=\bvi{x}{y}_{[x_P,x^P]}$ in case $x\equiv_Py$, and
 \begin{equation}\label{Eq:DefNewMeasQ}
 \begin{aligned}
 \bvQ{x}{y}=\bvP{x^P}{y_P}
 &\wedge
 \bigl(\bvP{x_P}{y_P}\vee\bve{x}{x_P}_{[x_P,x^P]}\bigr)\\
 &\wedge
 \bigl(\bvP{x^P}{y^P}\vee\bve{y^P}{y}_{[y_P,y^P]}\bigr)\\
 &\wedge
 \bigl(\bvP{x_P}{y^P}\vee\bve{x}{x_P}_{[x_P,x^P]}
 \vee\bve{y^P}{y}_{[y_P,y^P]}\bigr),
 \end{aligned}
 \end{equation}
if $x\not\equiv_Py$.

The proof of the following lemma is straightforward.

\begin{lemma}\label{L:bvQxy1inP}
The new map $\bvi{{}_{-}}{{}_{-}}$ extends
the original one $\bvi{{}_{-}}{{}_{-}}_{\bP}$, and also all maps
of the form $\bvi{{}_{-}}{{}_{-}}_{[x_P,x^P]}$, for $x\in Q$.
Furthermore, for all $x,y\in Q$, the following statements hold:
\begin{enumerate}
\item $x\in P$ implies that
$\bvQ{x}{y}=\bvP{x}{y_P}\wedge\bigl(
\bvP{x}{y^P}\vee\bve{y^P}{y}_{[y_P,y^P]}\bigr)$;

\item $y\in P$ implies that
$\bvQ{x}{y}=\bvP{x^P}{y}\wedge\bigl(
\bvP{x_P}{y}\vee\bve{x}{x_P}_{[x_P,x^P]}\bigr)$.
\end{enumerate}
\end{lemma}

\begin{lemma}\label{L:bvQleq=0}
$x\leq y$ implies that $\bvQ{x}{y}=0$, for all $x,y\in Q$.
\end{lemma}

\begin{proof}
If $x\equiv_Py$, then $\bvQ{x}{y}=\bvi{x}{y}_{[x_P,x^P]}=0$. If
$x\not\equiv_Py$, then, as $P\leint Q$, we get $x^P\leq y_P$,
thus $\bvP{x^P}{y_P}=0$, and so $\bvQ{x}{y}=0$.
\end{proof}

\begin{lemma}\label{L:1/3TRbvQ}
The inequality $\bvQ{x}{z}\leq\bvQ{x}{y}\vee\bvQ{y}{z}$ holds,
for all $x,y,z\in Q$ two of which belong to $P$.
\end{lemma}

\begin{proof}
Suppose first that $x,y\in P$. By applying
Lemma~\ref{L:bvQxy1inP} to $\bvQ{x}{z}$ and $\bvQ{y}{z}$, we
reduce the problem to the two inequalities
 \begin{align*}
 \bvP{x}{z_P}&\leq\bvP{x}{y}\wedge\bvP{y}{z_P}\,,\\
 \bvP{x}{z^P}&\leq\bvP{x}{y}\wedge\bvP{y}{z^P}\,,
 \end{align*}
which hold by assumption. The proof is dual for the case
$y,z\in P$.

Suppose now that $x,z\in P$. By applying
Lemma~\ref{L:bvQxy1inP} to $\bvQ{x}{y}$ and $\bvQ{y}{z}$, we
reduce the problem to four inequalities, which we proceed
to verify:
 \begin{align*}
 \bvP{x}{z}&\leq\bvP{x}{y_P}\vee\bvP{y_P}{z}\\
 &\leq\bvP{x}{y_P}\vee\bvP{y^P}{z}\,.
 \end{align*}
 \begin{align*}
 \bvP{x}{z}&\leq\bvP{x}{y_P}\vee\bvP{y_P}{z}\\
 &\leq\bvP{x}{y_P}\vee\bvP{y_P}{z}\vee\bve{y}{y_P}_{[y_P,y^P]}\,.
 \end{align*}
 \begin{align*}
 \bvP{x}{z}&\leq\bvP{x}{y^P}\vee\bvP{y^P}{z}\\
 &\leq\bvP{x}{y^P}\vee\bvP{y^P}{z}\vee\bve{y^P}{y}_{[y_P,y^P]}\,.
 \end{align*}
 \begin{align*}
 \bvP{x}{z}&
 \leq\bvP{x}{y^P}\vee\bve{y^P}{y_P}_{\bP}\vee\bvP{y_P}{z}\\
 &=\bvP{x}{y^P}\vee\bve{y^P}{y}_{[y_P,y^P]}
 \vee\bve{y}{y_P}_{[y_P,y^P]}\vee\bvP{y_P}{z}\,.
 \end{align*}
Hence $\bvP{x}{z}\leq\bvQ{x}{y}\vee\bvQ{y}{z}$\,.
\end{proof}

\begin{lemma}\label{L:UBxlePy}
The Boolean value $\bvQ{x}{y}$ lies below each of the
semilattice elements
$\bvP{x^P}{y_P}$\, ,
$\bvP{x_P}{y_P}\vee\bve{x}{x_P}_{[x_P,x^P]}$,
$\bvP{x^P}{y^P}\vee\bve{y^P}{y}_{[y_P,y^P]}$,
and $\bvP{x_P}{y^P}\vee\bve{x}{x_P}_{[x_P,x^P]}
\vee\bve{y^P}{y}_{[y_P,y^P]}$, for all $x,y\in Q$.
\end{lemma}

\begin{proof}
This is obvious by the definition of $\bvQ{x}{y}$ in case
$x\not\equiv_Py$. If $x\equiv_Py$, then, putting $u=x_P=y_P$ and
$v=x^P=y^P$, the four semilattice elements in the statement
above are respectively equal to $\bve{v}{u}_{\bP}$,
$\bve{x}{u}_{[u,v]}$, $\bve{v}{y}_{[u,v]}$, and
$\bve{x}{u}_{[u,v]}\vee\bve{v}{y}_{[u,v]}$. As
$\bve{v}{u}_{\bP}=\bve{v}{u}_{[u,v]}$, we need to prove that
$\bvi{x}{y}_{[u,v]}$ lies below both
$\bve{x}{u}_{[u,v]}$ and $\bve{v}{y}_{[u,v]}$, which is obvious
(for example,
$\bvi{x}{y}_{[u,v]}\leq\bvi{x}{u}_{[u,v]}\vee\bvi{u}{y}_{[u,v]}
=\bve{x}{u}_{[u,v]}$).
\end{proof}

\begin{lemma}\label{L:EncadrbvQxy}
The inequalities $\bvQ{x_P}{y}\leq\bvQ{x}{y}\leq\bvQ{x^P}{y}$
and\linebreak
$\bvQ{x}{y^P}\leq\bvQ{x}{y}\leq\bvQ{x}{y_P}$ hold, for all
$x,y\in Q$.
\end{lemma}

\begin{proof}
As the two sets of inequalities are dual, it suffices to
prove that\linebreak
$\bvQ{x_P}{y}\leq\bvQ{x}{y}\leq\bvQ{x^P}{y}$. As the conclusion
is obvious in case $x\equiv_Py$, it suffices to consider the
case where $x\not\equiv_Py$. As, by Lemma~\ref{L:bvQxy1inP}, the
equality
 \[
 \bvQ{x^P}{y}=\bvP{x^P}{y_P}\wedge\bigl(
 \bvP{x^P}{y^P}\vee\bve{y^P}{y}_{[y_P,y^P]}\bigr)
 \]
holds, it follows from Lemma~\ref{L:UBxlePy} that
$\bvQ{x}{y}\leq\bvQ{x^P}{y}$. Moreover, again by
Lemma~\ref{L:bvQxy1inP}, the equality
 \[
 \bvQ{x_P}{y}=\bvP{x_P}{y_P}\wedge\bigl(
 \bvP{x_P}{y^P}\vee\bve{y^P}{y}_{[y_P,y^P]}\bigr)
 \]
holds, and so $\bvQ{x_P}{y}$ lies below each of the four
meetands defining $\bvQ{x}{y}$ on the right hand side of
\eqref{Eq:DefNewMeasQ}, and hence $\bvQ{x_P}{y}\leq\bvQ{x}{y}$.
\end{proof}

\begin{lemma}\label{L:OtherEncbvQ}
The inequalities
$\bvQ{x}{y}\leq\bvQ{x_P}{y}\vee\bve{x}{x_P}_{[x_P,x^P]}$ and
$\bvQ{x}{y}\leq\bvQ{x}{y^P}\vee\bve{y^P}{y}_{[y_P,y^P]}$ hold,
for all $x,y\in Q$.
\end{lemma}

\begin{proof}
By symmetry, it suffices to prove the first inequality.
Using the expression of $\bvQ{x_P}{y}$ given by
Lemma~\ref{L:bvQxy1inP}(i), we reduce the problem to the
following two inequalities,
 \begin{align*}
 \bvQ{x}{y}&\leq\bvP{x_P}{y_P}\vee\bve{x}{x_P}_{[x_P,x^P]}\,,\\
 \bvQ{x}{y}&\leq\bvP{x_P}{y^P}\vee\bve{x}{x_P}_{[x_P,x^P]}
 \vee\bve{y^P}{y}_{[y_P,y^P]}\,,
 \end{align*}
that follow immediately from Lemma~\ref{L:UBxlePy}.
\end{proof}

\begin{lemma}\label{L:TwoMorebvQ}
The inequalities
$\bvQ{x^P}{y}\leq\bve{x^P}{x}_{[x_P,x^P]}\vee\bvQ{x}{y}$ and
$\bvQ{x}{y_P}\leq\bve{y}{y_P}_{[y_P,y^P]}\vee\bvQ{x}{y}$ hold,
for all $x,y\in Q$.
\end{lemma}

\begin{proof}
By symmetry, it suffices to prove the first inequality. Suppose
first that $x\equiv_Py$, put $u=x_P=y_P$ and $v=x^P=y^P$. We need
to prove that
$\bve{v}{y}_{\bQ}\leq\bve{v}{x}_{\bQ}\vee\bvi{x}{y}_{[u,v]}$,
which is obvious since $\bve{v}{y}_{\bQ}=\bve{v}{y}_{[u,v]}$ and
$\bve{v}{x}_{\bQ}=\bve{v}{x}_{[u,v]}$.

Now suppose that $x\not\equiv_Py$. As in \eqref{Eq:DefNewMeasQ},
$\bvQ{x}{y}$ is the meet of four meetands, so the first
inequality reduces to four inequalities, which we proceed to
prove:
 \begin{align*}
 \bvQ{x^P}{y}&\leq\bvP{x^P}{y_P}
 &&(\text{by Lemma~\ref{L:UBxlePy}})\\
 &\leq\bve{x^P}{x}_{\bQ}\vee\bvQ{x}{y_P}
 &&(\text{by Lemma~\ref{L:1/3TRbvQ}})\\
 &\leq\bve{x^P}{x}_{[x_P,x^P]}\vee\bvP{x^P}{y_P}
 &&(\text{by Lemma~\ref{L:EncadrbvQxy}}).
 \end{align*}
(We have used the easy observation that
$\bve{x^P}{x}_{\bQ}=\bve{x^P}{x}_{[x_P,x^P]}$.)
 \begin{align*}
 \bvQ{x^P}{y}&\leq\bvP{x^P}{y_P}
 \qquad(\text{by Lemma~\ref{L:EncadrbvQxy}})\\
 &\leq\bve{x^P}{x_P}_{\bP}\vee\bvP{x_P}{y_P}\\
 &=\bve{x^P}{x}_{[x_P,x^P]}\vee\bvP{x_P}{y_P}\vee
 \bve{x}{x_P}_{[x_P,x^P]}\,.
 \end{align*}
 \begin{align*}
 \bvQ{x^P}{y}&\leq\bvP{x^P}{y^P}\vee\bve{y^P}{y}_{\bQ}
 \qquad(\text{by Lemma~\ref{L:1/3TRbvQ}})\\
 &\leq\bve{x^P}{x}_{[x_P,x^P]}\vee\bvP{x^P}{y^P}\vee
 \bve{y^P}{y}_{[y_P,y^P]}\,.
 \end{align*}
 \begin{align*}
 \bvQ{x^P}{y}&\leq\bvP{x^P}{y^P}\vee\bve{y^P}{y}_{\bQ}
 \qquad(\text{by Lemma~\ref{L:1/3TRbvQ}})\\
 &\leq\bve{x^P}{x_P}_{\bP}\vee\bvP{x_P}{y^P}\vee
 \bve{y^P}{y}_{[y_P,y^P]}\\
 &=\bve{x^P}{x}_{[x_P,x^P]}\vee\bvP{x_P}{y^P}\vee
 \bve{x}{x_P}_{[x_P,x^P]}\vee\bve{y^P}{y}_{[y_P,y^P]}\,,
 \end{align*}
which completes the proof of the inequality
 \[
 \bvQ{x^P}{y}\leq\bve{x^P}{x}_{[x_P,x^P]}\vee\bvQ{x}{y}\,.
 \]
The proof of the inequality
$\bvQ{x}{y_P}\leq\bve{y}{y_P}_{[y_P,y^P]}\vee\bvQ{x}{y}$ is dual.
\end{proof}

\begin{lemma}\label{L:TrIneqbvQ}
The inequality $\bvQ{x}{z}\leq\bvQ{x}{y}\vee\bvQ{y}{z}$ holds,
for all $x,y,z\in Q$.
\end{lemma}

\begin{proof}
This is obvious in case $x\equiv_Py\equiv_Pz$, as
$\bvi{{}_{-}}{{}_{-}}_{[x_P,x^P]}$ is a p-measure. So
suppose that either $x\not\equiv_Py$ or $y\not\equiv_Pz$, say
$x\not\equiv_Py$. Expressing the Boolean value $\bvQ{x}{y}$ as
in \eqref{Eq:DefNewMeasQ}, we reduce the problem to four
inequalities, that we proceed to prove:
 \begin{align*}
 \bvQ{x}{z}&\leq\bvQ{x^P}{z}
 &&(\text{by Lemma~\ref{L:EncadrbvQxy}})\\
 &\leq\bvP{x^P}{y_P}\vee\bvQ{y_P}{z}
 &&(\text{by Lemma~\ref{L:1/3TRbvQ}})\\
 &\leq\bvP{x^P}{y_P}\vee\bvQ{y}{z}
 &&(\text{by Lemma~\ref{L:EncadrbvQxy}}).
 \end{align*}
 \begin{align*}
 \bvQ{x}{z}&\leq\bvQ{x_P}{z}\vee\bve{x}{x_P}_{[x_P,x^P]}
 &&(\text{by Lemma~\ref{L:OtherEncbvQ}})\\
 &\leq\bvP{x_P}{y_P}\vee\bvQ{y_P}{z}\vee\bve{x}{x_P}_{[x_P,x^P]}
 &&(\text{by Lemma~\ref{L:1/3TRbvQ}})\\
 &\leq\bvP{x_P}{y_P}\vee\bve{x}{x_P}_{[x_P,x^P]}\vee\bvQ{y}{z}
 &&(\text{by Lemma~\ref{L:EncadrbvQxy}}).
 \end{align*}
 \begin{align*}
 \bvQ{x}{z}&\leq\bvQ{x^P}{z}
 &&(\text{by Lemma~\ref{L:EncadrbvQxy}})\\
 &\leq\bvP{x^P}{y^P}\vee\bvQ{y^P}{z}
 &&(\text{by Lemma~\ref{L:1/3TRbvQ}})\\
 &\leq\bvP{x^P}{y^P}\vee\bve{y^P}{y}_{[y_P,y^P]}
 \vee\bvQ{y}{z}
 &&(\text{by Lemma~\ref{L:TwoMorebvQ}}).
 \end{align*}
 \begin{align*}
 \bvQ{x}{z}&\leq\bvQ{x_P}{z}\vee\bve{x}{x_P}_{[x_P,x^P]}
 \hskip3cm(\text{by Lemma~\ref{L:OtherEncbvQ}})\\
 &\leq\bvP{x_P}{y^P}\vee\bvP{y^P}{z}\vee\bve{x}{x_P}_{[x_P,x^P]}
 \qquad(\text{by Lemma~\ref{L:1/3TRbvQ}})\\
 &\leq\bvP{x_P}{y^P}\vee\bve{y^P}{y}_{[y_P,y^P]}
 \vee\bve{x}{x_P}_{[x_P,x^P]}\vee\bvQ{y}{z}\\
 &\hskip8cm
 (\text{by Lemma~\ref{L:TwoMorebvQ}}).
 \end{align*}
This completes the proof.
\end{proof}

So we have reached the following result.

\begin{proposition}\label{P:1ExtMeas}
Let $D$ be a distributive lattice with zero, let $\bP$ be a
$D$-valued p-measured poset, and let $Q$ be an
interval extension of $P$ in which each
interval of~$Q$ of the form $[x_P,x^P]$, for $x\in Q$, is
endowed with a p-measure $\bvi{{}_{-}}{{}_{-}}_{[x_P,x^P]}$
such that
$\bve{x^P}{x_P}_{\bP}=\bve{x^P}{x_P}_{[x_P,x^P]}$ for all
$x\in Q$. Then there exists a common extension of all p-measures
$\bvi{{}_{-}}{{}_{-}}_{\bP}$ and
$\bvi{{}_{-}}{{}_{-}}_{[x_P,x^P]}$, for $x\in Q$, to a
p-measure on $Q$, given by \eqref{Eq:DefNewMeasQ} on pairs
$\seq{x,y}$ such that $x\not\equiv_Py$.
\end{proposition}

\section{Doubling extensions; the conditions (DB1) and (DB2)}
\label{S:Doubling}

For a poset $\Lambda$ and a $\Lambda$-indexed diagram
$\vec{S}=\seqm{S_i,\varphi_{i,j}}{i\leq j\text{ in }\Lambda}$ of
\jzs s and \jzh s, we shall say that a $\Lambda$-indexed diagram
$\seqm{\bQ_i}{i\in\Lambda}$ of p-measured posets is
\emph{$\vec{S}$-valued}, if $\bQ_i$ is $S_i$-valued and
$\bQ_j$ extends $\bQ_i$ with respect to~$\varphi_{i,j}$ for all
$i\leq j$ in $\Lambda$.

We shall also use the convention of notation and
terminology that consists of extending to p-measured
posets the notions defined for posets, by restricting them to
the underlying posets and stating that the poset extensions
involved preserve the corresponding p-measures.
For example, we say that a p-measured poset~$\bQ$ is an
\emph{interval extension} of a p-measured poset~$\bP$,
in notation $\bP\leint\bQ$, if~$\bQ$ extends~$\bP$ and the
underlying posets (cf. Notation~\ref{Not:BVs}) satisfy
$P\leint Q$. In particular, a normal interval diagram of
p-measured lattices is a diagram of p-measured lattices whose
underlying posets form a normal interval diagram.

\begin{definition}\label{D:Doubling}
Let $\bP$ and $\bQ$ be p-measured posets such that
$\bP\lerc\bQ$. We say that~$\bQ$ is a \emph{doubling
extension} of~$\bP$, in notation $\bP\ledb\bQ$,
if $\bve{x}{x_P}\sim\bve{x^P}{x}$ for all
$x\in Q$. Equivalently, either
$\bve{x^P}{x_P}=\bve{x}{x_P}$ or $\bve{x^P}{x_P}=\bve{x^P}{x}$,
for all $x\in Q$.
\end{definition}

The following lemma shows that under mild assumptions, doubling
extensions are transitive.

\begin{lemma}\label{L:ExtDoubl}
Let $\bP$, $\bQ$, and $\bR$ be p-measured
posets. If $P\lerc Q\lerc R$, $P\leint R$, and
$\bP\ledb\bQ\ledb\bR$, then $\bP\ledb\bR$.
\end{lemma}

\begin{proof}
Let $x\in R$, we prove $\bve{x}{x_P}\sim\bve{x^P}{x}$. As
$\bQ\ledb\bR$, we get $\bve{x}{x_Q}\sim\bve{x^Q}{x}$. Hence, if
$\set{x_Q,x^Q}\subseteq P$, then $x_P=x_Q$ and $x^P=x^Q$, thus
we are done. Suppose that $\set{x_Q,x^Q}\not\subseteq P$, say
$x^Q\notin P$. As $P\leint R$ and $x\sim x^Q$, it follows from
Lemma~\ref{L:compxxP} that $x\sim(x^Q)_P$. If $x\leq(x^Q)_P$,
then, as $(x^Q)_P\leq x^Q$ and $(x^Q)_P$ belongs to~$P$, we get
$x^Q=(x^Q)_P\in P$, \contr; hence $(x^Q)_P\leq x$. As
$x_P\leq(x^Q)_P$ and $(x^Q)_P\in P$, we get $(x^Q)_P=x_P$. As
$\bP\ledb\bQ$, we get $\bve{x^P}{x^Q}\sim\bve{x^Q}{x_P}$. But
this also holds trivially in case $x^Q\in P$, so it holds in
every case. So we have proved the following:
 \begin{equation}\label{Eq:1bveComp}
 \bve{x^P}{x^Q}\sim\bve{x^Q}{x_P}.
 \end{equation}
The dual argument gives
 \begin{equation}\label{Eq:2bveComp}
 \bve{x_Q}{x_P}\sim\bve{x^P}{x_Q}.
 \end{equation}
If $\bve{x^P}{x_Q}\leq\bve{x_Q}{x_P}$, then we get
$\bve{x_Q}{x_P}=\bve{x^P}{x_P}$, and thus
$\bve{x}{x_P}=\bve{x^P}{x_P}$, and we are done. Dually, the same
conclusion follows from $\bve{x^Q}{x_P}\leq\bve{x^P}{x^Q}$.

By \eqref{Eq:1bveComp} and \eqref{Eq:2bveComp}, it remains to
consider the case where both inequalities
$\bve{x_Q}{x_P}\leq\bve{x^P}{x_Q}$ and
$\bve{x^P}{x^Q}\leq\bve{x^Q}{x_P}$ hold, in which case
 \begin{equation}\label{Eq:dDdbv-s}
 \bve{x^Q}{x_P}=\bve{x^P}{x_Q}=\bve{x^P}{x_P}.
 \end{equation}
{}From $\bQ\ledb\bR$ it follows that
$\bve{x^Q}{x}\sim\bve{x}{x_Q}$. Suppose, for example, that
$\bve{x}{x_Q}\leq\bve{x^Q}{x}$. Hence
$\bve{x^Q}{x}=\bve{x^Q}{x_Q}$, and we get
 \begin{align*}
 \bve{x^P}{x}&=\bve{x^P}{x^Q}\vee\bve{x^Q}{x}\\
 &=\bve{x^P}{x^Q}\vee\bve{x^Q}{x_Q}\\
 &=\bve{x^P}{x_Q}\\
 &=\bve{x^P}{x_P}&&(\text{see \eqref{Eq:dDdbv-s}})\\
 &\geq\bve{x}{x_P}.\tag*{\qed}
 \end{align*}
\renewcommand{\qed}{}
\end{proof}

{}From now on until the end of this section, we shall
fix a finite lattice~$\Lambda$ with largest element~$\ell$, a
$\Lambda$-indexed diagram
$\vec{D}=\seqm{D_i,\varphi_{i,j}}{i\leq j\text{ in }\Lambda}$ of
distributive lattices with zero and \jzh s, a
$\vec{D}\res_{<\ell}$\,-valued normal interval diagram
$\seqm{\bQ_i}{i<\ell}$ of p-measured lattices. In addition, we
assume that the following statements hold:
\begin{itemize}
\item[(DB1)] $\bQ_j$ is a doubling extension of
$\bQ_i$ for all $i\leq j<\ell$.

\item[(DB2)] For all $i<\ell$ and all $x,y\in Q_i$ with
$\nu(x)\nleq\nu(y)$,
$\bve{x}{x_{\bl}}_{\bQ_i}=\bve{x^{\bl}}{x_{\bl}}_{\bQ_i}$
implies that $\bvi{x}{y}_{\bQ_i}=\bvi{x^{\bl}}{y}_{\bQ_i}$ and
$\bve{x^{\bl}}{x}_{\bQ_i}=\bve{x^{\bl}}{x_{\bl}}_{\bQ_i}$
implies that $\bvi{y}{x}_{\bQ_i}=\bvi{y}{x_{\bl}}_{\bQ_i}$.
\end{itemize}

As usual, we denote by $P=\bigcup\famm{Q_i}{i<\ell}$ the strong
amalgam of $\seqm{Q_i}{i<\ell}$.

\begin{remark}\label{Rk:Add3inc}
It suffices to verify (DB2) in case $x\inc y$. Indeed, let
$x,y\in Q_i$ such that $\nu(x)\nleq\nu(y)$
and $\bve{x}{x_{\bl}}=\bve{x^{\bl}}{x_{\bl}}$
(to ease the notation, we drop the indices~$\bQ_i$). If $x\leq
y$, then, by Lemma~\ref{L:Existx**}, $x^{\bl}\leq y$, thus
$\bvi{x}{y}=\bvi{x^{\bl}}{y}=0$. If $y\leq x$, then, again by
Lemma~\ref{L:Existx**}, $y\leq x_{\bl}$, and so
 \begin{align*}
 \bvi{x}{y}&=\bve{x}{x_{\bl}}\vee\bve{x_{\bl}}{y}
 &&(\text{because }y\leq x_{\bl}\leq x)\\
 &=\bve{x^{\bl}}{x_{\bl}}\vee\bve{x_{\bl}}{y}
 &&(\text{because }\bve{x^{\bl}}{x_{\bl}}=\bve{x}{x_{\bl}})\\
 &=\bvi{x^{\bl}}{y}.
 \end{align*}
The proof that $x\sim y$ and $\nu(x)\nleq\nu(y)$ and
$\bve{x^{\bl}}{x}=\bve{x^{\bl}}{x_{\bl}}$
implies that $\bvi{y}{x}=\bvi{y}{x_{\bl}}$ is dual.
\end{remark}

\begin{notation}\label{Not:Poplusmin}
We add a largest element, denoted by $1$, to $D_{\ell}$,
and for all $x,y\in P$, we define an element
$\bvo{x}{y}$ of $D_{\ell}\cup\set{1}$ as follows:
 \begin{equation}\label{Eq:Defbvoxy}
 \bvo{x}{y}=\begin{cases}
 \varphi_{i,\ell}\bigl(\bvi{x}{y}_{\bQ_i}\bigr),&
 \text{ if }\nu(x)\vee\nu(y)\leq i<\ell,\\
 1,&\text{otherwise}.
 \end{cases}
 \end{equation}
It is obvious that the value of $\bvo{x}{y}$ defined in the
first case is independent of the choice of~$i$ such that
$\nu(x)\vee\nu(y)\leq i<\ell$. We also put
 \[
 \bveo{x}{y}=\bvo{\max\set{x,y}}{\min\set{x,y}},
 \qquad\text{for all comparable }x,y\in P.
 \]

\begin{lemma}\label{L:PdoublesQi}
The elements $\bveo{x}{\lc{x}{i}}$ and $\bveo{\uc{x}{i}}{x}$ are
comparable, for all $x\in P$ and all $i<\ell$. Furthermore,
$\bveo{x^{\bl}}{x}\sim\bveo{x}{x_{\bl}}$ for all
$x\in P\setminus Q_0$.
\end{lemma}

\begin{proof}
Let $x\in P$. As $x\in Q_j$ for some $j<\ell$, we get
$\lc{x}{i}=\lc{x}{i\wedge j}$ and $\uc{x}{i}=\uc{x}{i\wedge j}$
(cf. Lemma~\ref{L:PresRCINT}(i)). As $\bQ_{i\wedge j}\ledb\bQ_j$,
we get $\bve{x}{\lc{x}{i\wedge j}}_{\bQ_j}\sim
\bve{\uc{x}{i\wedge j}}{x}_{\bQ_j}$, that is,
$\bve{x}{\lc{x}{i}}_{\bQ_j}\sim\bve{\uc{x}{i}}{x}_{\bQ_j}$, and
thus, applying $\varphi_{j,\ell}$, we obtain the relation
$\bveo{x}{\lc{x}{i}}\sim\bveo{\uc{x}{i}}{x}$.

It follows from Lemma~\ref{L:Existx**} that $\nu(x_{\bl})$ and
$\nu(x^{\bl})$ are comparable and that, if~$i$ denotes their
maximum, then $x_{\bl}=\lc{x}{i}$ and $x^{\bl}=\uc{x}{i}$. By
applying the result of the previous paragraph, we obtain
$\bveo{x^{\bl}}{x}\sim\bveo{x}{x_{\bl}}$.
\end{proof}

Now we put
 \begin{align*}
 P^{\oplus}&=\setm{x\in P\setminus Q_0}
 {\bveo{x}{x_{\bl}}=\bveo{x^{\bl}}{x_{\bl}}},\\
 P^{\ominus}&=\setm{x\in P\setminus Q_0}
 {\bveo{x^{\bl}}{x}=\bveo{x^{\bl}}{x_{\bl}}}.
 \end{align*}
\end{notation}
If $x$ belongs to $Q_i\setminus Q_0$, then
$x^{\bl},x_{\bl}\in Q_i$. Hence, both $\bveo{x}{x_{\bl}}$ and
$\bveo{x^{\bl}}{x}$ are evaluated by the formula giving the case
$\nu(x)\vee\nu(y)<\ell$ of \eqref{Eq:Defbvoxy}. Therefore, it
follows from Lemma~\ref{L:PdoublesQi} that
$P\setminus Q_0=P^{\oplus}\cup P^{\ominus}$.

\section{Strong amalgams of p-measured posets; from $\bvo{x}{y}$
to $\bvi{x}{y}$}\label{S:variousbv}

{}From now on until Lemma~\ref{L:PropagADD3}, we shall
fix a finite lattice~$\Lambda$ with largest element~$\ell$, a
$\Lambda$-indexed diagram
$\vec{D}=\seqm{D_i,\varphi_{i,j}}{i\leq j\text{ in }\Lambda}$ of
\emph{finite} distributive lattices and \jzh s, a
$\vec{D}\res_{<\ell}$\,-valued normal interval diagram
$\seqm{\bQ_i}{i<\ell}$ of p-measured lattices.
Furthermore, we assume that the conditions (DB1) and (DB2)
introduced in Section~\ref{S:Doubling} are satisfied.

We denote by $P$ the strong amalgam of
$\seqm{Q_i}{i<\ell}$ and by $\rho(x)$ the height of~$\nu(x)$ in
$\Lambda$, for all $x\in P$.

\begin{lemma}\label{L:FirstTR3bv}
For every positive integer $n$ and all elements
$x_0,x_1,\dots,x_n\in P$,
$\nu(x_0)\vee\nu(x_n)<\ell$ implies that
$\bvo{x_0}{x_n}\leq\bigvee_{i<n}\bvo{x_i}{x_{i+1}}$.
\end{lemma}

\begin{proof}
We argue by induction on the pair
$\seq{n,\sum_{k=0}^n\rho(x_k)}$, ordered lexicographically.
The conclusion is trivial for $n=1$.

Now suppose that $n=2$. If either $\nu(x_0)\vee\nu(x_1)=\ell$ or
$\nu(x_1)\vee\nu(x_2)=\ell$, then the right hand side of the
desired inequality is equal to~$1$ and we are done; so suppose
that $\nu(x_0)\vee\nu(x_1),\nu(x_1)\vee\nu(x_2)<\ell$.
If $\nu(x_1)\leq\nu(x_0)\vee\nu(x_2)$, then, putting\linebreak
$k=\nu(x_0)\vee\nu(x_2)$ (which is smaller than~$\ell$),
all the Boolean values under consideration are images under
$\varphi_{k,\ell}$ of the corresponding Boolean values
in~$\bQ_k$, so the conclusion follows from the inequality
$\bvi{x_0}{x_2}_{\bQ_k}\leq
\bvi{x_0}{x_1}_{\bQ_k}\vee\bvi{x_1}{x_2}_{\bQ_k}$
(we will often encounter this kind of reduction, and
we will summarize it by ``everything happens below level~$k$'').
Now suppose that $\nu(x_1)\nleq\nu(x_0)\vee\nu(x_2)$. In
particular, $x_1\notin Q_0$ and $\nu(x_1)\nleq\nu(x_0),\nu(x_2)$.
By Lemma~\ref{L:PdoublesQi}, $x_1$ belongs to
$P^{\oplus}\cup P^{\ominus}$. If $x_1\in P^{\oplus}$, then, as
$\nu(x_1)\vee\nu(x_2)<\ell$ and by (DB2),
$\bvo{x_1}{x_2}=\bvo{(x_1)^{\bl}}{x_2}$, hence
 \begin{align*}
 \bvo{x_0}{x_2}&\leq
 \bvo{x_0}{(x_1)^{\bl}}\vee\bvo{(x_1)^{\bl}}{x_2}
 &&(\text{by the induction hypothesis})\\
 &\leq\bvo{x_0}{x_1}\vee\bvo{x_1}{x_2}
 &&(\text{because }\bvo{x_0}{(x_1)^{\bl}}\leq\bvo{x_0}{x_1})
 \end{align*}
so we are done. The proof is symmetric in case
$x_1\in P^{\ominus}$. This concludes the case where $n=2$.

Now assume that $n\geq3$. It $\nu(x_i)\vee\nu(x_{i+1})=\ell$ for
some $i<n$, then the right hand side of the desired inequality is
equal to~$1$ and we are done; so suppose that
$\nu(x_i)\vee\nu(x_{i+1})<\ell$ for all $i<n$. Suppose that
there are $i,j$ such that $0\leq i\leq j\leq n$ and $2\leq j-i<n$
such that $\nu(x_i)\vee\nu(x_j)<\ell$. It follows from the
induction hypothesis that
$\bvo{x_i}{x_j}\leq\bigvee_{i\leq k<j}\bvo{x_k}{x_{k+1}}$. Hence,
using again the induction hypothesis, we get
 \begin{align*}
 \bvo{x_0}{x_n}&\leq\bigvee_{k<i}\bvo{x_k}{x_{k+1}}\vee
 \bvo{x_i}{x_j}\vee\bigvee_{j\leq k<n}\bvo{x_k}{x_{k+1}}\\
 &\leq\bigvee_{k<n}\bvo{x_k}{x_{k+1}}\,,
 \end{align*}
so we are done again. Hence suppose that
$0\leq i\leq j\leq n$ and $2\leq j-i<n$ implies that
$\nu(x_i)\vee\nu(x_j)=\ell$, for all $i,j$. As
$\nu(x_1)\vee\nu(x_n)=\ell$ while $\nu(x_0)\vee\nu(x_n)<\ell$
(\emph{we use here the assumption that $n\geq3$}),
we get $\nu(x_1)\nleq\nu(x_0)$. As
$\nu(x_2)\vee\nu(x_3)<\ell$ and $\nu(x_1)\vee\nu(x_3)=\ell$, we
get $\nu(x_1)\nleq\nu(x_2)$. Hence, if $x_1\in P^{\oplus}$, then,
as $\nu(x_1)\nleq\nu(x_2)$ and by (DB2),
$\bvo{x_1}{x_2}=\bvo{(x_1)^{\bl}}{x_2}$, and hence, by using the
induction hypothesis and the obvious inequality
$\bvo{x_0}{(x_1)^{\bl}}\leq\bvo{x_0}{x_1}$
(``everything there happens below level
$\nu(x_0)\vee\nu(x_1)$''), we get
 \begin{align*}
 \bvo{x_0}{x_n}&\leq\bvo{x_0}{(x_1)^{\bl}}\vee
 \bvo{(x_1)^{\bl}}{x_2}\vee
 \bigvee_{2\leq i<n}\bvo{x_i}{x_{i+1}}\\
 &\leq\bigvee_{k<n}\bvo{x_k}{x_{k+1}}\,,
 \end{align*}
so we are done.
If $x_1\in P^{\ominus}$, then, as $\nu(x_1)\nleq\nu(x_0)$ and by
(DB2), $\bvo{x_0}{x_1}=\bvo{x_0}{(x_1)_{\bl}}$, hence, by
using the induction hypothesis and the obvious inequality
$\bvo{(x_1)_{\bl}}{x_2}\leq\bvo{x_1}{x_2}$, we get
 \begin{align*}
 \bvo{x_0}{x_n}&\leq\bvo{x_0}{(x_1)_{\bl}}\vee
 \bvo{(x_1)_{\bl}}{x_2}\vee
 \bigvee_{2\leq i<n}\bvo{x_i}{x_{i+1}}\\
 &\leq\bigvee_{k<n}\bvo{x_k}{x_{k+1}}\,,
 \end{align*}
so we are done. As $x_1\in P^{\oplus}\cup P^{\ominus}$ (cf.
Lemma~\ref{L:PdoublesQi}), this completes the induction step.
\end{proof}

\begin{notation}\label{Not:variousBV}
We put
 \begin{gather*}
 P(z)=\setm{t\in P}{\nu(t)<\nu(z)},\\
 P^{\oplus}(z)=P(z)\cap P^{\oplus},\qquad
 P^{\ominus}(z)=P(z)\cap P^{\ominus},
 \end{gather*}
for all $z\in P$.
Furthermore, for all $x,y\in P$, we define
 \begin{align}
 \bvp{x}{y}&=\bigwedge\Famm
 {\bvo{x}{t}\vee\bvo{t}{y}}
 {t\in P(y)}\,,\label{Eq:Defbvp}\\
 \bvm{x}{y}&=\bigwedge\Famm
 {\bvo{x}{t}\vee\bvo{t}{y}}
 {t\in P(x)}\,,\label{Eq:Defbvm}\\
 \bvpm{x}{y}&=\bigwedge\Famm
 {\bvo{x}{u}\vee\bvo{u}{v}\vee\bvo{v}{y}}
 {\seq{u,v}\in P^{\ominus}(x)\times P^{\oplus}(y)}\,,
 \label{Eq:Defbvpm}\\
 \bvi{x}{y}
 &=\bvo{x}{y}\wedge\bvp{x}{y}\wedge
 \bvm{x}{y}\wedge\bvpm{x}{y}\,.\label{Eq:Defbvd}
 \end{align}
(All meets are evaluated in $D_{\ell}$, the empty meet
being defined as equal to~$1$.)
We observe that the meet on the right hand side of
\eqref{Eq:Defbvp} may be taken over all $t\in P(y)$
such that $\nu(x)\vee\nu(t)<\ell$: indeed, for all other
$t\in P(y)$, we get $\bvo{x}{t}=1$. Similarly, the
meet on the right hand side of \eqref{Eq:Defbvm} may be taken
over all $t\in P(x)$ such that
$\nu(t)\vee\nu(y)<\ell$, and the meet on the right hand side of
\eqref{Eq:Defbvpm} may be taken over all
$\seq{u,v}\in P^{\ominus}(x)\times P^{\oplus}(y)$ such that
$\nu(u)\vee\nu(v)<\ell$.
\end{notation}

\begin{lemma}\label{L:FirstTrbvi}
$\bvi{x}{y}\leq\bvo{x}{t}\vee\bvo{t}{y}$ for all $x,y,t\in P$.
\end{lemma}

\begin{proof}
We argue by induction on $\rho(x)+\rho(y)+\rho(t)$. If
$\nu(x)\vee\nu(t)=\ell$ or\linebreak $\nu(t)\vee\nu(y)=\ell$ then
the right hand side of the desired inequality is equal to~$1$ so
we are done. Suppose, from now on, that
$\nu(x)\vee\nu(t),\nu(t)\vee\nu(y)<\ell$. If
$\nu(x)\vee\nu(y)<\ell$, then it follows from
Lemma~\ref{L:FirstTR3bv} (for $n=2$) that
$\bvo{x}{y}\leq\bvo{x}{t}\vee\bvo{t}{y}$, so we are done as
$\bvi{x}{y}\leq\bvo{x}{y}$. Now suppose that
$\nu(x)\vee\nu(y)=\ell$. In particular, $x,y\notin Q_0$. If
$t\in Q_0$, then $t\in P(y)$, thus
 \[
 \bvi{x}{y}\leq\bvp{x}{y}\leq\bvo{x}{t}\vee\bvo{t}{y}\,.
 \]
So suppose that $t\notin Q_0$. If $\nu(x)\leq\nu(t)$, then
``everything happens below level $\nu(t)\vee\nu(y)$'' (which is
smaller than~$\ell$), so we are done. The conclusion is similar
in case $\nu(y)\leq\nu(t)$.

So suppose that $\nu(x),\nu(y)\nleq\nu(t)$. If
$\nu(t)\leq\nu(y)$, then $\nu(t)<\nu(y)$ (because
$\nu(y)\nleq\nu(t)$), thus $t\in P(y)$, and thus
 \[
 \bvi{x}{y}\leq\bvp{x}{y}\leq\bvo{x}{t}\vee\bvo{t}{y}\,,
 \]
so we are done. If $t\in P^{\oplus}$ and $\nu(t)\nleq\nu(y)$,
then, by (DB2), $\bvo{t}{y}=\bvo{t^{\bl}}{y}$, and thus
 \begin{align*}
 \bvi{x}{y}&\leq\bvo{x}{t^{\bl}}\vee\bvo{t^{\bl}}{y}
 &&(\text{by the induction hypothesis})\\
 &\leq\bvo{x}{t}\vee\bvo{t}{y}\,,
 \end{align*}
so we are done again. This covers the case where
$t\in P^{\oplus}$. The proof is symmetric for $t\in P^{\ominus}$.
\end{proof}

\begin{lemma}\label{L:2ndTrbvi}
$\bvi{x}{y}\leq\bvo{x}{u}\vee\bvo{u}{v}\vee\bvo{v}{y}$ for all
$x,y,u,v\in P$.
\end{lemma}

\begin{proof}
We argue by induction on $\rho(x)+\rho(y)+\rho(u)+\rho(v)$.
If either $\nu(x)\vee\nu(u)=\ell$ or $\nu(u)\vee\nu(v)=\ell$ or
$\nu(v)\vee\nu(y)=\ell$, then the right hand side of the desired
inequality is equal to~$1$ and we are done. So suppose that
$\nu(x)\vee\nu(u),\nu(u)\vee\nu(v),\nu(v)\vee\nu(y)<\ell$. If
$\nu(x)\vee\nu(v)<\ell$, then, by Lemma~\ref{L:FirstTR3bv}, we
get $\bvo{x}{v}\leq\bvo{x}{u}\vee\bvo{u}{v}$, and so
 \begin{align*}
 \bvi{x}{y}&\leq\bvo{x}{v}\vee\bvo{v}{y}
 &&(\text{by Lemma~\ref{L:FirstTrbvi}})\\
 &\leq\bvo{x}{u}\vee\bvo{u}{v}\vee\bvo{v}{y}\,.
 \end{align*}
The conclusion is similar for $\nu(u)\vee\nu(y)<\ell$. So
suppose that $\nu(x)\vee\nu(v)=\nu(u)\vee\nu(y)=\ell$. In
particular, $\nu(x)\nleq\nu(u)$,
$\nu(y)\nleq\nu(v)$, and $\nu(u)\inc\nu(v)$, so
$x,y,u,v\notin Q_0$.

Suppose that $u\in P^{\oplus}$. As $\nu(u)\nleq\nu(v)$ and by
(DB2), we get $\bvo{u}{v}=\bvo{u^{\bl}}{v}$, hence
 \begin{align*}
 \bvi{x}{y}&\leq
 \bvo{x}{u^{\bl}}\vee\bvo{u^{\bl}}{v}\vee\bvo{v}{y}
 &&(\text{by the induction hypothesis})\\
 &\leq\bvo{x}{u}\vee\bvo{u}{v}\vee\bvo{v}{y}
 &&(\text{because }\bvo{x}{u^{\bl}}\leq\bvo{x}{u}).
 \end{align*}
Suppose that $u\in P^{\ominus}$ and $\nu(u)\not<\nu(x)$. As
$\nu(x)\nleq\nu(u)$, we get $\nu(u)\nleq\nu(x)$, thus
$\bvo{x}{u}=\bvo{x}{u_{\bl}}$, and so
 \begin{align*}
 \bvi{x}{y}&\leq
 \bvo{x}{u_{\bl}}\vee\bvo{u_{\bl}}{v}\vee\bvo{v}{y}
 &&(\text{by the induction hypothesis})\\
 &\leq\bvo{x}{u}\vee\bvo{u}{v}\vee\bvo{v}{y}
 &&(\text{because }\bvo{u_{\bl}}{v}\leq\bvo{u}{v}).
 \end{align*}
The case where either $v\in P^{\ominus}$ or ($v\in P^{\oplus}$
and $\nu(v)\not<\nu(y)$ is symmetric. The only remaining case is
where $u\in P^{\ominus}(x)$ and $v\in P^{\oplus}(y)$, in which
case
 \begin{equation*}
 \bvi{x}{y}\leq\bvpm{x}{y}\leq
 \bvo{x}{u}\vee\bvo{u}{v}\vee\bvo{v}{y}.\tag*{\qed}
 \end{equation*}
\renewcommand{\qed}{}
\end{proof}

Consequently, we get the formula
 \begin{equation}\label{Eq:ShorterDefbvi}
 \bvi{x}{y}=\bigwedge
 \Famm{\bvo{x}{u}\vee\bvo{u}{v}\vee\bvo{v}{y}}{u,v\in P}\,,
 \qquad\text{for all }x,y\in P.
 \end{equation}

\begin{lemma}\label{L:HalfTR4bvd}
$\bvi{x}{z}\leq\bvi{x}{y}\vee\bvo{y}{z}$ for all $x,y,z\in P$.
\end{lemma}

\begin{proof}
If $\nu(y)\vee\nu(z)=\ell$ then $\bvo{y}{z}=1$ and the
conclusion is trivial. Suppose that $\nu(y)\vee\nu(z)<\ell$. A
direct use of Lemma~\ref{L:FirstTrbvi} yields the inequality
$\bvi{x}{z}\leq\bvo{x}{y}\vee\bvo{y}{z}$, while a direct use of
Lemma~\ref{L:2ndTrbvi} together with the distributivity
of~$D_{\ell}$ yields that
$\bvi{x}{z}\leq\bvp{x}{y}\wedge\bvm{x}{y}\vee\bvo{y}{z}$. It
remains to establish the inequality
$\bvi{x}{z}\leq\bvpm{x}{y}\vee\bvo{y}{z}$, which reduces, by the
distributivity of~$D_{\ell}$, to proving the inequality
 \begin{equation}\label{Eq:ixzoxuvyz}
 \bvi{x}{z}\leq
 \bvo{x}{u}\vee\bvo{u}{v}\vee\bvo{v}{y}\vee\bvo{y}{z},
 \end{equation}
for all $\seq{u,v}\in P^{\ominus}(x)\times P^{\oplus}(y)$. {}From
$\nu(v)<\nu(y)$ it follows that
$\bvo{v}{z}\leq\bvo{v}{y}\vee\bvo{y}{z}$ (``everything
there happens below level $\nu(y)\vee\nu(z)$''), and hence
 \begin{align*}
 \bvi{x}{z}&\leq\bvo{x}{u}\vee\bvo{u}{v}\vee\bvo{v}{z}
 &&(\text{by Lemma~\ref{L:2ndTrbvi}})\\
 &\leq\bvo{x}{u}\vee\bvo{u}{v}\vee\bvo{v}{y}\vee\bvo{y}{z},
 \end{align*}
which completes the proof of \eqref{Eq:ixzoxuvyz}.
\end{proof}

\begin{lemma}\label{L:FullTRbvd}
$\bvi{x}{z}\leq\bvi{x}{y}\vee\bvi{y}{z}$ for all $x,y,z\in P$.
\end{lemma}

\begin{proof}
For elements $u,v\in P$, we get, by three successive applications
of Lemma~\ref{L:HalfTR4bvd}, the inequalities
 \begin{align*}
 \bvi{x}{u}&\leq\bvi{x}{y}\vee\bvo{y}{u}\,;\\
 \bvi{x}{v}&\leq\bvi{x}{u}\vee\bvo{u}{v}\,;\\
 \bvi{x}{z}&\leq\bvi{x}{v}\vee\bvo{v}{z}\,.
 \end{align*}
Hence, combining these inequalities, we obtain
 \[
 \bvi{x}{z}\leq
 \bvi{x}{y}\vee\bvo{y}{u}\vee\bvo{u}{v}\vee\bvo{v}{z}\,.
 \]
Evaluating the meets of both sides over $u,v\in P$
and using (the easy direction of) \eqref{Eq:ShorterDefbvi}
yields the desired conclusion.
\end{proof}

As a consequence, we obtain the following simple expression of
$\bvi{x}{y}$.

\begin{corollary}\label{C:Simplebvd}
The Boolean value $\bvi{x}{y}$ is equal to the meet in
$D_{\ell}$ of all elements of $D_{\ell}$ of the form
 \begin{equation}\label{Eq:joinbvzizi+1}
 \bvo{x}{z_1}\vee\bvo{z_1}{z_2}\vee\cdots\vee\bvo{z_{n-1}}{y}\,,
 \end{equation}
where $n$ is a natural number and $z_0,z_1,\dots,z_n\in P$
such that $z_0=x$, $z_n=y$, and $\nu(z_i)\vee\nu(z_{i+1})<\ell$
for all $i<n$. Furthermore, it is sufficient to restrict the
meet to finite sequences $\seq{z_0,z_1,z_2,z_3}$ \pup{so $n=3$}.
\end{corollary}

\begin{proof}
Denote temporarily by $\bvi{x}{y}^*$ the meet in $D_{\ell}$ of
all elements of $D_{\ell}$ of the form \eqref{Eq:joinbvzizi+1}.
An immediate application of the easy direction of
\eqref{Eq:ShorterDefbvi} yields the inequality
$\bvi{x}{y}^*\leq\bvi{x}{y}$. Conversely, for every natural
number $n$ and all $z_0,z_1,\dots,z_n\in P$
such that $z_0=x$, $z_n=y$, and $\nu(z_i)\vee\nu(z_{i+1})<\ell$
for all $i<n$,
 \begin{align*}
 \bvi{x}{y}&\leq\bigvee_{i<n}\bvi{z_i}{z_{i+1}}
 &&(\text{by Lemma~\ref{L:FullTRbvd}})\\
 &\leq\bigvee_{i<n}\bvo{z_i}{z_{i+1}}
 &&(\text{because }\bvi{z_i}{z_{i+1}}\leq\bvo{z_i}{z_{i+1}}),
 \end{align*}
which concludes the proof of the first part. The bound $n=3$
follows from the easy direction of
\eqref{Eq:ShorterDefbvi}.
\end{proof}

As an immediate consequence of Lemma~\ref{L:FirstTR3bv}, we
obtain that the equality\linebreak $\bvi{x}{y}=\bvo{x}{y}$ holds
for all $x,y\in P$ such that $\nu(x)\vee\nu(y)<\ell$. Hence we
obtain the following lemma.

\begin{lemma}\label{L:bvdextbvi}
The p-measure $\bvi{{}_{-}}{{}_{-}}$ extends the p-measure
$\bvi{{}_{-}}{{}_{-}}_{\bQ_i}$ with respect to
$\varphi_{i,\ell}$, for all $i<\ell$.
\end{lemma}

\begin{definition}\label{D:StrongAmalgMP}
The strong amalgam $P=\bigcup\famm{Q_i}{i<\ell}$, endowed
with the p-measure $\bvi{{}_{-}}{{}_{-}}$ constructed above,
will be called the \emph{strong amalgam of the family
$\seqm{\bQ_i}{i<\ell}$ with respect to~$\vec{D}$}.
\end{definition}

So we have reached the main goal of the present section.

\begin{proposition}\label{P:PresTR3Amalg}
Let $\Lambda$ be a finite lattice with largest element~$\ell$,
let\linebreak
$\vec{D}=\seqm{D_i,\varphi_{i,j}}{i\leq j\text{ in }\Lambda}$ be a
$\Lambda$-indexed diagram of finite distributive lattices and
\jzh s, and let $\seqm{\bQ_i}{i<\ell}$ be a
$\vec{D}\res_{<\ell}$\,-valued normal interval diagram of p-measured
lattices satisfying
\textup{(DB1)} and \textup{(DB2)}. Then the strong amalgam
$\bP$ of $\seqm{\bQ_i}{i<\ell}$ \pup{see
Definition~\textup{\ref{D:StrongAmalgMP}}} is a $D_{\ell}$-valued
p-measured lattice, which extends $\bQ_i$ with respect to
$\varphi_{i,\ell}$, for all $i<\ell$.
\end{proposition}

\begin{lemma}\label{L:bPdbQi}
Under the assumptions of
Proposition~\textup{\ref{P:PresTR3Amalg}}, the p-measured
poset~$\bP$ is a doubling extension of~$\bQ_i$ for all $i<\ell$.
\end{lemma}

\begin{proof}
An immediate consequence of Lemmas~\ref{L:PdoublesQi}
and~\ref{L:bvdextbvi}.
\end{proof}

The goal of the following lemma is to propagate the
assumption~(DB2) through the induction process that will appear
in the constructions of Theorems~\ref{T:ConstrMeasDiagr}
and~\ref{T:PosetRepr}.

\begin{lemma}\label{L:PropagADD3}
For all $x,y\in P$, the following statements hold:
\begin{enumerate}
\item \pup{$\nu(x)\nleq\nu(y)$ and $x\in P^{\oplus}$} implies
that $\bvi{x}{y}=\bvi{x^{\bl}}{y}$.

\item \pup{$\nu(y)\nleq\nu(x)$ and $y\in P^{\ominus}$} implies
that $\bvi{x}{y}=\bvi{x}{y_{\bl}}$.
\end{enumerate}
\end{lemma}

\begin{proof}
As (i) and (ii) are dual, it suffices to establish (i).
We first claim that for all $x\in P^{\oplus}$ and all $y\in P$,
$\nu(x)\nleq\nu(y)$ implies that
$\bvo{x^{\bl}}{y}\leq\bvo{x}{y}$. Indeed, the equality holds by
assumption (DB2) in case $\nu(x)\vee\nu(y)<\ell$. If
$\nu(x)\vee\nu(y)=\ell$, then $\bvo{x}{y}=1$ and we are done
again.

Now let $x\in P^{\oplus}$ and $y\in P$ such that
$\nu(x)\nleq\nu(y)$, we must prove that $\bvi{x^{\bl}}{y}$ lies
below $\bvo{x}{y}$, $\bvp{x}{y}$, $\bvm{x}{y}$, and
$\bvpm{x}{y}$.
 \begin{align*}
 \bvi{x^{\bl}}{y}&\leq\bvo{x^{\bl}}{y}
 &&(\text{see \eqref{Eq:Defbvd}})\\
 &\leq\bvo{x}{y}
 &&(\text{as }\nu(x)\nleq\nu(y)
 \text{ and by the claim above}). \end{align*}
Now let $t\in P(y)$.
 \begin{align*}
 \bvi{x^{\bl}}{y}&\leq\bvo{x^{\bl}}{t}\vee\bvo{t}{y}
 &&(\text{by Lemma~\ref{L:FirstTrbvi}})\\
 &\leq\bvo{x}{t}\vee\bvo{t}{y}
 &&(\text{as }\nu(x)\nleq\nu(t)
 \text{ and by the claim above}).
 \end{align*}
Evaluating the meets of both sides over $t\in P^{\oplus}(y)$
yields $\bvi{x^{\bl}}{y}\leq\bvp{x}{y}$. A similar (but not
symmetric!) proof yields $\bvi{x^{\bl}}{y}\leq\bvm{x}{y}$.
Finally, let
$u\in P^{\ominus}(x)$ and $v\in P^{\oplus}(y)$. Then
 \begin{align*}
 \bvi{x^{\bl}}{y}&\leq
 \bvo{x^{\bl}}{u}\vee\bvo{u}{v}\vee\bvo{v}{y}
 &&(\text{by Lemma~\ref{L:2ndTrbvi}})\\
 &\leq\bvo{x}{u}\vee\bvo{u}{v}\vee\bvo{v}{y}
 &&(\text{as }\nu(x)\nleq\nu(u)
 \text{ and by the claim above}),
 \end{align*}
hence, evaluating the meets of both sides over
$\seq{u,v}\in P^{\ominus}(x)\times P^{\oplus}(y)$, we get
$\bvi{x^{\bl}}{y}\leq\bvpm{x}{y}$, which completes the proof.
\end{proof}

\section{Constructing a p-measure on a covering, doubling
extension of a strong amalgam of lattices}\label{S:MeasQell}

Let $\Lambda$, $\vec{D}$, and $\seqm{\bQ_i}{i<\ell}$ satisfy the
assumptions of Proposition~\ref{P:PresTR3Amalg}, with
strong amalgam~$\bP$ (cf. Definition~\ref{D:StrongAmalgMP}).
By Proposition~\ref{P:PresTR3Amalg}, $\bvP{{}_{-}}{{}_{-}}$ is a
$D_{\ell}$-valued p-measure on~$\bP$, which extends
each p-measured lattice $\bQ_i$ with respect to the corresponding
\jzh\ $\varphi_{i,\ell}$.

Now we let $Q$ be a covering extension of $P$.
Furthermore, we assume that each closed interval $[x_P,x^P]$ of
$Q$, for $x\in Q$, is endowed with a p-measure
$\bvi{{}_{-}}{{}_{-}}_{[x_P,x^P]}$ such that
 \begin{align}
 \bve{x}{x_P}_{[x_P,x^P]}&\sim\bve{x^P}{x}_{[x_P,x^P]}\,,
 &&\text{for all }x\in Q,
 \label{Eq:Complevelell}\\
 \bve{x^P}{x_P}_{[x_P,x^P]}&=\bve{x^P}{x_P}_{\bP}\,,
 &&\text{for all }x\in Q.
 \label{Eq:Cohell}
 \end{align}
(Observe that the notation $\bve{x^P}{x_P}_{\bP}$ in
\eqref{Eq:Cohell} above does not involve
the full definition of the strong amalgam given in
Definition~\ref{D:StrongAmalgMP}: indeed, from
$P\lecov Q$ it follows that $x_P\preceq_Px^P$; as $P$ is the
strong amalgam of $\seqm{Q_i}{i<\ell}$, $x_P$ and $x^P$ belong to
some $Q_i$, and so we can just put
$\bve{x^P}{x_P}_{\bP}=
\varphi_{i,\ell}\bigl(\bve{x^P}{x_P}_{\bQ_i}\bigr)$,
which is independent of the chosen $i$.)

The goal of the present section is to extend
$\bvP{{}_{-}}{{}_{-}}$ to a p-measure on~$\bQ$ such
that, setting $\bQ_{\ell}=\bQ$, the extended diagram
$\seqm{\bQ_i}{i\leq\ell}$ satisfies the assumptions of
Proposition~\ref{P:PresTR3Amalg}.

We need to verify several points. First, for all $i<\ell$, as
$Q_i\leint P$ and $P\lecov Q$, we obtain from
Lemma~\ref{L:PresRCINT} that $Q_i\leint Q$. Item~(3) of
Definition~\ref{D:SAsys} for the extended
diagram $\seqm{\bQ_i}{i\leq\ell}$ follows from the definition of
the ordering of~$P_{\ell}$ (cf. Section~\ref{S:SASIS}).
Further, the new valuation on the extended diagram
$\seqm{\bQ_i}{i\leq\ell}$ extends the original one (so we shall
still denote it by~$\nu$), and $\nu(x)=\ell$ for all
$x\in Q\setminus P$. In addition, the elements $x_{\bl}$ and
$x^{\bl}$ (cf. Lemma~\ref{L:Existx**}) remain the same for
$x\in P\setminus Q_0$, while $x_{\bl}=x_P$ and $x^{\bl}=x^P$ for
all $x\in Q\setminus P$.

Now we denote by $\bvQ{{}_{-}}{{}_{-}}$ the p-measure
that we constructed in Section~\ref{S:ExtMeasInt} (cf.
Proposition~\ref{P:1ExtMeas}), extending $\bvP{{}_{-}}{{}_{-}}$
and all p-measures
$\bvi{{}_{-}}{{}_{-}}_{[x_P,x^P]}$, for $x\in Q$---this is made
possible by \eqref{Eq:Cohell}. It follows from the assumption
\eqref{Eq:Complevelell} that
$\bve{x^P}{x}_{\bQ}\sim\bve{x}{x_P}_{\bQ}$ for all $x\in Q$;
that is, \emph{$\bQ$ is a doubling extension of~$\bP$}. As
$Q_i\leint Q$ and by Lemmas~\ref{L:bPdbQi}
and~\ref{L:ExtDoubl} (applied to the extensions
$\bQ_i\leq\bP\leq\bQ$), we obtain that \emph{$\bQ$ is a doubling
extension of $\bQ_i$}. This takes care of extending (DB1) to
the larger diagram.

It remains to verify that $\seqm{\bQ_i}{i\leq\ell}$ satisfies
(DB2). So let $x,y\in Q$ such that $\nu(x)\nleq\nu(y)$, we need
to verify that
$\bve{x}{x_{\bl}}_{\bQ}=\bve{x^{\bl}}{x_{\bl}}_{\bQ}$ implies
that $\bvQ{x}{y}=\bvQ{x^{\bl}}{y}$ and
$\bve{x^{\bl}}{x}_{\bQ}=\bve{x^{\bl}}{x_{\bl}}_{\bQ}$
implies that $\bvQ{y}{x}=\bvQ{y}{x_{\bl}}$. We prove for example
the first statement. {}From $\nu(x)\nleq\nu(y)$ it follows that
$y\in P$. If $x\in P$ then we are done by
Lemma~\ref{L:PropagADD3}, so the remaining case is where
$x\in Q\setminus P$. Observe that $x_{\bl}=x_P$ and
$x^{\bl}=x^P$. As $y\in P$, the Boolean value $\bvQ{x}{y}$ is
given by Lemma~\ref{L:bvQxy1inP}(ii). Hence proving the
inequality $\bvQ{x^{\bl}}{y}\leq\bvQ{x}{y}$ reduces to proving
that $\bvQ{x^P}{y}$ (of course equal to $\bvP{x^P}{y}$)
lies below both $\bvP{x^P}{y}$ and
$\bvP{x_P}{y}\vee\bve{x}{x_P}_{[x_P,x^P]}$. The first inequality
is a tautology, and the second one is proved as follows:
 \begin{align*}
 \bvP{x^P}{y}&\leq\bve{x^P}{x_P}_{\bP}\vee\bvP{x_P}{y}
 &&(\text{because }\bvP{{}_{-}}{{}_{-}}\text{ is a p-measure})\\
 &=\bve{x}{x_P}_{\bQ}\vee\bvP{x_P}{y}
 &&(\text{because }
 \bve{x}{x_{\bl}}_{\bQ}=\bve{x^{\bl}}{x_{\bl}}_{\bQ})\\
 &=\bve{x}{x_P}_{[x_P,x^P]}\vee\bvP{x_P}{y}\,.
 \end{align*}
As the inequality $\bvQ{x}{y}\leq\bvQ{x^{\bl}}{y}$ always holds,
we have proved the equality, and hence
the extended diagram $\seqm{\bQ_i}{i\leq\ell}$ satisfies (DB2).
So we have reached the following theorem, which is the main
technical result of the present paper. It refers to the
conditions (DB1) and (DB2) introduced in
Section~\ref{S:Doubling}.

\begin{theorem}\label{T:ConstrMeasDiagr}
Let $\Lambda$ be a finite lattice with largest element~$\ell$,
let\linebreak
$\vec{D}=\seqm{D_i,\varphi_{i,j}}{i\leq j\text{ in }\Lambda}$ be a
$\Lambda$-indexed diagram of finite distributive lattices and
\jzh s, and let $\seqm{\bQ_i}{i<\ell}$ be a
$\vec{D}\res_{<\ell}$\,-valued normal interval diagram of p-measured
lattices satisfying
\textup{(DB1)} and \textup{(DB2)}. Let~$Q$ be a covering
extension of the strong amalgam~$P$ of
$\seqm{Q_i}{i<\ell}$. Furthermore, we assume that for all
$x\in Q$, the closed interval $[x_P,x^P]$ of $Q$ is endowed with
a p-measure $\bvi{{}_{-}}{{}_{-}}_{[x_P,x^P]}$ 
\pup{depending only of the interval $[x_P,x^P]$} such that
 \[
 \bve{x}{x_P}_{[x_P,x^P]}\sim\bve{x^P}{x}_{[x_P,x^P]}
 \quad\text{and}\quad\bve{x^P}{x_P}_{[x_P,x^P]}=
 \bve{x^P}{x_P}_{\bP}\,.
 \]
Then there exists a
$D_{\ell}$-valued p-measure on~$Q$ extending all p-measures
$\bvi{{}_{-}}{{}_{-}}_{[x_P,x^P]}$ such that,
defining $\bQ_{\ell}$ as the corresponding p-measured poset, the
extended diagram $\seqm{\bQ_i}{i\leq\ell}$ is a
$\vec{D}$-valued normal interval diagram of p-measured posets
satisfying \textup{(DB1)} and \textup{(DB2)}.
\end{theorem}

This result makes it possible to state and prove our main
theorem.

\begin{theorem}\label{T:PosetRepr}
Let $\Lambda$ be a lower finite \ms\ and let\linebreak
$\vec{D}=\seqm{D_i,\varphi_{i,j}}{i\leq j\text{ in }\Lambda}$ be a
$\Lambda$-indexed diagram of finite distributive lattices and
\jzuh s. Then there exists a $\vec{D}$-valued normal
interval diagram $\seqm{\bQ_i}{i\in\Lambda}$ of finite p-measured
lattices satisfying \textup{(DB1)} and \textup{(DB2)} together
with the following additional conditions:
\begin{enumerate}
\item For all $i<j$ in $\Lambda$ and all $x<y$ in $Q_i$, there
exists $z\in Q_j$ such that $x<z<y$.

\item $\bve{y}{x}_{\bQ_i}\in\J(D_i)\cup\set{0}$, for
all $i\in\Lambda$ and all $x,y\in Q_i$ such that $x\prec_{Q_i}y$.

\item For all $i\in\Lambda$ and all $p\in\J(D_i)$, there exists
$x\in Q_i$ such that $0\prec_{Q_i}x$ and $\bve{x}{0}_{\bQ_i}=p$.
\end{enumerate}
\end{theorem}

\begin{proof}
We construct $\bQ_i$ by induction on the height of $i$
in $\Lambda$. After possibly adding a new zero element to
$\Lambda$, we may assume that $D_0=\set{0,1}$, so we take
$Q_0=\set{0,1}$, with the p-measure defined by
$\bve{1}{0}_{\bQ_0}=1$. Put
$\Lambda_n=\setm{i\in\Lambda}{\hgt(i)\leq n}$ and denote by
$\vec{D}_{(n)}$ the restriction of~$\vec{D}$ to $\Lambda_n$, for every
natural number~$n$. Suppose having constructed a
$\vec{D}_{(n)}$-valued normal interval diagram
$\seqm{\bQ_i}{i\in\Lambda_n}$ of finite p-measured lattices
satisfying (DB1), (DB2), and Conditions (i)--(iii) of the
statement of the theorem, we show how to extend it to a
$\vec{D}_{(n+1)}$-valued normal interval diagram of finite p-measured
lattices satisfying (DB1) and (DB2). In order to propagate
Item~(2) of Definition~\ref{D:SAsys} through our induction,
we shall add the following induction hypothesis:
\begin{equation}\label{Eq:Qxyldisj}
\text{Every }x\in\bigcup\famm{Q_i}{i\in\Lambda_n}
\text{ can be written in the form }\seq{\ol{x},\nu(x)},
\end{equation}
where $\nu$ denotes the valuation associated with the diagram
$\seqm{Q_i}{i\in\Lambda_n}$.
Let $\ell\in\Lambda_{n+1}\setminus\Lambda_n$ and denote
by~$\bP_{\ell}$ the strong amalgam of $\seqm{\bQ_i}{i<\ell}$
with respect to~$\vec{D}$ given in Definition~\ref{D:StrongAmalgMP}.
It follows from Proposition~\ref{P:SA2forP} that $P_{\ell}$ is a
lattice and every~$Q_i$, for $i<\ell$, is a sublattice
of~$P_{\ell}$. For all $x\prec y$ in $P_{\ell}$ such that $\bve{y}{x}_{\bP_\ell}>0$, we put
 \begin{align*}
 1_{x,y,\ell}&
 =\setm{p\in\J(D_{\ell})}{p\leq\bve{y}{x}_{\bP_{\ell}}},\\
 B_{x,y,\ell}&=\Pow(1_{x,y,\ell})
 \qquad\qquad(\text{the powerset lattice of }1_{x,y,\ell}),\\
 \ol{Q}_{x,y,\ell}&=\setm{\seq{X,\es}}{X\subseteq 1_{x,y,\ell}}
 \cup\setm{\seq{1_{x,y,\ell},Y}}{Y\subseteq 1_{x,y,\ell}}.
 \end{align*}
Observe that the condition $\bve{y}{x}_{\bP_\ell}>0$ implies that the set $1_{x,y,\ell}$, that we shall often denote by~$1$, is nonempty. Also, $\ol{Q}_{x,y,\ell}$ is a sublattice of $B_{x,y,\ell}\times B_{x,y,\ell}$.

Hence $\ol{Q}_{x,y,\ell}$ is the ordinal sum of two
copies of the Boolean lattice~$B_{x,y,\ell}$, with the top of the
lower copy of $B_{x,y,\ell}$ (namely, $\seq{X,\es}$ where
$X=1_{x,y,\ell}$) identified with the bottom of the upper copy of
$B_{x,y,\ell}$ (namely, $\seq{1,Y}$ where $Y=\es$).

We endow $\ol{Q}_{x,y,\ell}$ with the p-measure
$\bvi{{}_{-}}{{}_{-}}_{x,y,\ell}$ defined by
 \begin{align*}
 \bvi{\seq{X_0,\es}}{\seq{X_1,\es}}_{x,y,\ell}
 &=\bigvee(X_0\setminus X_1)\,,\\
 \bvi{\seq{1,Y_0}}
 {\seq{1,Y_1}}_{x,y,\ell}
 &=\bigvee(Y_0\setminus Y_1)\,,\\
 \bvi{\seq{X,\es}}{\seq{1,Y}}_{x,y,\ell}&=0\,,\\
 \bvi{\seq{1,Y}}{\seq{X,\es}}_{x,y,\ell}&=
 \bigvee\bigl(\complement X\cup Y\bigr),
 \end{align*}
(where we put $\complement X=1_{x,y,\ell}\setminus X$),
for all $X,X_0,X_1,Y,Y_0,Y_1\subseteq 1$ (it is easy to verify
that this way we get, indeed, a p-measure on
$\ol{Q}_{x,y,\ell}$). Further, we put
 \begin{align*}
 Q'_{x,y,\ell}&=\ol{Q}_{x,y,\ell}\setminus
 \set{\seq{\es,\es},\seq{1,1}}
 &&(\text{`truncated }\ol{Q}_{x,y,\ell}\text{'}),\\
 Q_{x,y,\ell}&=\setm{\seq{\seq{t,x,y},\ell}}{t\in Q'_{x,y,\ell}},
 \end{align*}
where $Q'_{x,y,\ell}$ is endowed with the restrictions of both
the ordering and the p-measure of $\ol{Q}_{x,y,\ell}$ and
$Q_{x,y,\ell}$ is endowed with the ordering and p-measure for
which the map $t\mapsto\seq{\seq{t,x,y},\ell}$ is a measure-preserving
isomorphism. So $Q_{x,y,\ell}$ is the result of applying to
$\ol{Q}_{x,y,\ell}$ the following two transformations:
\begin{itemize}
\item[---] Remove the top and bottom elements of
$\ol{Q}_{x,y,\ell}$; get $Q'_{x,y,\ell}$.

\item[---] Replace $t$ by $\seq{\seq{t,x,y},\ell}$, for all
$t\in Q'_{x,y,\ell}$; get $Q_{x,y,\ell}$.
\end{itemize}
The latter step (from $Q'_{x,y,\ell}$ to~$Q_{x,y,\ell}$) is put there in order to ensure the induction hypothesis~\eqref{Eq:Qxyldisj} while making the~$Q_{x,y,\ell}$s pairwise disjoint.

In case $\bve{y}{x}_{\bP_\ell}=0$, we pick an outside element~$t_{x,y,\ell}$ and we set $Q_{x,y,\ell}=\set{\seq{t_{x,y,\ell},\ell}}$, the one-element poset. Furthermore, we endow $\ol{Q}_{x,y,\ell}=\set{x,y,\seq{t_{x,y,\ell},\ell}}$ with the p-measure with constant value~$0$.

Observe that in any case, $Q_{x,y,\ell}$ is nonempty.

Put $Q_{\ell}=P_{\ell}+\sum\famm{Q_{x,y,\ell}}
{x\prec y\text{ in }P_{\ell}}$
(cf. \eqref{Eq:SumNotSISch}). Then~\eqref{Eq:Qxyldisj} is maintained at level~$\ell$, and $Q_{\ell}$ is an interval
extension of~$P_{\ell}$ (cf. Lemma~\ref{L:SumSISch}). In fact, as~$Q_{x,y,\ell}$ is defined only for $x\prec y$ in~$P_{\ell}$,
the poset~$Q_{\ell}$ is a covering extension of~$P_{\ell}$ (cf.
Definition~\ref{D:covering}). We shall still denote by
$\bvi{{}_{-}}{{}_{-}}_{x,y,\ell}$ the p-measure on
$Q_{x,y,\ell}\cup\set{x,y}$ induced by the p-measure on
$\ol{Q}_{x,y,\ell}$ defined above. As $P_{\ell}$ is a lattice and
$[x,y]_{Q_{\ell}}=Q_{x,y,\ell}\cup\set{x,y}\cong\ol{Q}_{x,y,\ell}$
is a lattice for all $x\prec y$ in $P_{\ell}$, it follows from
Lemma~\ref{L:IntLatt} that $Q_{\ell}$ is a lattice.

Now we verify Conditions \eqref{Eq:Complevelell} and
\eqref{Eq:Cohell} with respect to
$\bvi{{}_{-}}{{}_{-}}_{\bP_{\ell}}$ and all p-measures
$\bvi{{}_{-}}{{}_{-}}_{x,y,\ell}$. Fix $x\prec y$ in $P_{\ell}$ and let $z\in Q_{x,y,\ell}$; so $z_{P_\ell}=x$ and $z^{P_\ell}=y$. If $\bve{y}{x}_{\bP_\ell}=0$, then all members of both~\eqref{Eq:Complevelell} and~\eqref{Eq:Cohell} are zero, thus trivializing the corresponding statements. Hence suppose that $\bve{y}{x}_{\bP_\ell}>0$.
Condition~\eqref{Eq:Complevelell} follows immediately from the
inequalities
 \begin{align*}
 \bve{\seq{X,\es}}{\seq{\es,\es}}_{x,y,\ell}&=\bigvee X\leq
 \bve{y}{x}_{\bP_{\ell}}=\bve{\seq{1,1}}{\seq{X,\es}}_{x,y,\ell}\,,\\
 \bve{\seq{1,1}}{\seq{1,X}}_{x,y,\ell}&=
 \bigvee\bigl(\complement X\bigr)\leq\bve{y}{x}_{\bP_{\ell}}
 =\bve{\seq{1,X}}{\seq{\es,\es}}_{x,y,\ell}\,,
 \end{align*}
for all $X\subseteq 1_{x,y,\ell}$.
Condition \eqref{Eq:Cohell} follows from the equalities
 \[
 \bve{\seq{1,1}}{\seq{\es,\es}}_{x,y,\ell}=\bigvee 1_{x,y,\ell}
 =\bve{y}{x}_{\bP_{\ell}}\,.
 \]
Hence, by Theorem~\ref{T:ConstrMeasDiagr}, there is a
p-measure on $Q_{\ell}$, extending all p-measures
$\bvi{{}_{-}}{{}_{-}}_{x,y,\ell}$, such that
$\seqm{\bQ_i}{i\leq\ell}$ is a $\vec{D}_{\leq\ell}$\,-valued normal
interval diagram of p-measured lattices satisfying (DB1) and
(DB2).

Now we verify Conditions~(i)--(iii) of the statement of
Theorem~\ref{T:PosetRepr}. Let $i<\ell$ and let $x<y$ in $Q_i$,
we prove that $x\not\prec_{Q_{\ell}}y$. If
$x\not\prec_{P_{\ell}}y$ then this is trivial, so suppose that
$x\prec_{P_{\ell}}y$. Pick any element $z\in Q_{x,y,\ell}$ (we have seen that~$Q_{x,y,\ell}$ is always nonempty); then
$x<z<y$ in~$Q_{\ell}$. Condition~(i) follows.

In order to verify Condition~(ii) at level $Q_{\ell}$, it
suffices to prove that $\bve{v}{u}_{x,y,\ell}$
belongs to $\J(D_{\ell})\cup\set{0}$, for all $x\prec y$ in $P_{\ell}$ and
all $u\prec v$ in $\ol{Q}_{x,y,\ell}$. This is trivial in case $\bve{y}{x}_{\bP_\ell}=0$, in which case $\bve{v}{u}_{x,y,\ell}=0$. So suppose that $\bve{y}{x}_{\bP_\ell}>0$. There are a proper
subset $X$ of $1_{x,y,\ell}$ and an element $p\in\complement X$
such that either ($u=\seq{X,\es}$ and $v=\seq{X\cup\set{p},\es}$) or
($u=\seq{1,X}$ and $v=\seq{1,X\cup\set{p}}$). In both cases,
$\bve{v}{u}_{x,y,\ell}=p$ belongs to $\J(D_{\ell})$.

Now we verify Condition~(iii). Let $p\in\J(D_{\ell})$ and pick
$k\prec\ell$ in $\Lambda$. As
$p\leq\varphi_{k,\ell}(1)=
\bigvee\famm{\varphi_{k,\ell}(q)}{q\in\J(D_k)}$ and $p$ is
\jirr, there exists $q\in\J(D_k)$ such that
$p\leq\varphi_{k,\ell}(q)$. By the induction hypothesis
(Condition~(iii)), there exists $x\in Q_k$ such that
$0\prec_{Q_k}x$ and $\bve{x}{0}_{Q_k}=q$. Suppose that there
exists $y\in P_{\ell}$ such that $0<y<x$, and let $i<\ell$ such
that $y\in Q_i$. As $y\leq x$, there exists $z\in Q_{i\wedge k}$
such that $y\leq z\leq x$. As $0<z\leq x$ with $z\in Q_k$ and
$0\prec_{Q_k}x$, we get $z=x$, and so $x\in Q_{i\wedge k}$. If
$i\wedge k<k$, then, by Condition~(i) on $\vec{D}_{(n)}$, we get
$0\not\prec_{Q_k}x$, \contr. Therefore, $k=i\wedge k\leq i$, but
$k\prec\ell$, and thus $i=k$. As $y\in Q_k$, $0<y<x$, and
$0\prec_{Q_k}x$, we get again \contr. So we have proved that
$0\prec_{P_{\ell}}x$. As
$p\leq\varphi_{k,\ell}(q)=\bve{x}{0}_{\bP_{\ell}}$, we get
$p\in 1_{0,x,\ell}$. We consider the element
$t=\seq{\seq{\seq{\set{p},\es},0,x},\ell}$ of $Q_{0,x,\ell}$ (so
$0\prec t<x$ in $Q_{\ell}$). We compute
 \[
 \bve{t}{0}_{\bQ_{\ell}}=
 \bve{\seq{\set{p},\es}}{\seq{\es,\es}}_{0,x,\ell}=p\,,
 \]
which completes the verification of Condition~(iii) at
level~$\ell$.

In order to verify that $\seqm{\bQ_i}{i\in\Lambda_{n+1}}$ is as
required, it remains to verify
that $\seqm{Q_i}{i\in\Lambda_{n+1}}$ satisfies
Item~(2) of Definition~\ref{D:SAsys}. So let
$i,j\in\Lambda_{n+1}$, we need to verify that
$Q_i\cap Q_j=Q_{i\wedge j}$. This holds by induction hypothesis
for $i,j\in\Lambda_n$. As it trivially holds for $i=j$, we
assume that $i\neq j$. If $\hgt(i)=\hgt(j)=n$, then
 \begin{align}
 Q_i&=P_i\cup\bigcup\famm{Q_{x,y,i}}{x\prec y\text{ in }P_i}\,,
 \label{Eq:Qi}\\
 Q_j&=P_j\cup\bigcup\famm{Q_{x,y,j}}{x\prec y\text{ in }P_j}\,,
 \label{Eq:Qj}
 \end{align}
and thus, as $i\inc j$ and as~\eqref{Eq:Qxyldisj} is valid at all levels below either~$i$ or~$j$,
 \[
 Q_i\cap Q_j=P_i\cap P_j=
 \bigcup\famm{Q_{i'}\cap Q_{j'}}{i'<i,\ j'<j}=Q_{i\wedge j}\,.
 \]
If $\hgt(i)=n$ while $\hgt(j)<n$, then $Q_i$ is still given by
\eqref{Eq:Qi}, and so
 \[
 Q_i\cap Q_j=P_i\cap Q_j=\bigcup\famm{Q_{i'}\cap Q_j}{i'<i}
 =Q_{i\wedge j}\,,
 \]
which completes the verification of Item~(2) of
Definition~\ref{D:SAsys}. This completes the proof of
the induction step.
\end{proof}

\begin{remark}\label{Rk:PosetRepr}
More can be said in case all transition homomorphisms~$\varphi_{i,j}$ \emph{separate zero}, that is, $\varphi_{i,j}^{-1}\set{0}=\set{0}$, for all $i\leq j$ in~$\Lambda$. Indeed, in such a case, in the proof of Theorem~\ref{T:PosetRepr}, for all $x\prec y$ in~$P_\ell$, there exists $i<\ell$ such that $x,y\in Q_i$, and so $\bve{y}{x}_{\bP_\ell}=\varphi_{i,\ell}(\bve{y}{x}_{\bQ_i})$ is nonzero in case we have included in the induction hypothesis the assumption that $\bve{v}{u}>0$ for all~$i<\ell$ and all~$u<v$ in~$Q_i$. Hence, $\bve{y}{x}_{\bP_\ell}>0$ for all~$x<y$ in~$P_\ell$. Therefore, we can strengthen the conclusion~(ii) of Theorem~\ref{T:PosetRepr} by stating that $\bve{y}{x}$ is \jirr\ in~$D_i$, for all~$i\in\Lambda$ and all $x\prec y$ in~$Q_i$.
\end{remark}

\begin{corollary}\label{C:PosetRepr}
For every distributive \jzs\ $S$, there are a \mzs\ $P$ and a
$S$-valued p-measure $\bvi{{}_{-}}{{}_{-}}$ on
$P$ satisfying the following additional conditions:
\begin{enumerate}
\item $\bve{y}{x}>0$ for all $x<y$ in $P$.

\item For all $x\leq y$ in $P$ and all
$\ba,\bb\in S$, if $\bve{y}{x}\leq\ba\vee\bb$, there are a
positive integer~$n$ and a decomposition
$x=z_0\leq z_1\leq\dots\leq z_n=y$ such that either
$\bve{z_{i+1}}{z_i}\leq\ba$ or $\bve{z_{i+1}}{z_i}\leq\bb$, for
all $i<n$.

\item The subset $\setm{\bve{x}{0}}{x\in P}$ generates the
semilattice $S$.
\end{enumerate}
Furthermore, if $S$ is bounded, then $P$ can be taken
a bounded lattice.
\end{corollary}

\begin{proof}
Suppose first that $S$ is bounded.
By Lemma~\ref{L:ErsPudl}, $S$ is the directed union of its
finite distributive \jzu-subsemilattices. Hence we can write~$S$ as a directed union $S=\bigcup\famm{D_i}{i\in I}$, where $\Lambda$ is the
(lower finite) lattice of all finite subsets of~$S$,
all the $D_i$ are finite distributive $0,1$-subsemilattices of
$S$, and the transition map from $D_i$ to $D_j$ is the inclusion
map, for all $i\leq j$ in $\Lambda$ (in particular, it separates zero). Let
$\vec{\bQ}=\seqm{\bQ_i}{i\in\Lambda}$ be as in Theorem~\ref{T:PosetRepr}.
We prove that the union of all the p-measures
$\bvi{{}_{-}}{{}_{-}}_{\bQ_i}$ on
$Q=\bigcup\famm{Q_i}{i\in\Lambda}$ is as required.
Condition~(i) above follows from Remark~\ref{Rk:PosetRepr}. For Condition~(ii), suppose that $\bve{y}{x}\leq\ba\vee\bb$. Let~$i\in\Lambda$ such that $x,y\in Q_i$ and $\ba,\bb\in D_i$. As~$Q_i$ is finite, there exists a chain in~$Q_i$ of the form
 \[
 x=z_0\prec_{Q_i}z_1\prec_{Q_i}\cdots\prec_{Q_i}z_n=y.
 \]
For each~$i<n$, $\bve{z_{i+1}}{z_i}\leq\bve{y}{x}\leq\ba\vee\bb$. As~$\vec{\bQ}$ satisfies the conclusion of Theorem~\ref{T:PosetRepr}(ii), $\bve{z_{i+1}}{z_i}$ belongs to $\J(D_i)\cup\set{0}$, hence, as~$D_i$ is distributive, either $\bve{z_{i+1}}{z_i}\leq\ba$ or $\bve{z_{i+1}}{z_i}\leq\bb$. Condition~(ii) above follows.
As~$\J(D_i)$ join-generates~$D_i$, for each~$i\in\Lambda$, Condition~(iii) above follows from Condition~(iii) in Theorem~\ref{T:PosetRepr}.

In the general case, we apply the result above to
$S\cup\set{1}$ (for some new unit element $1$), and then,
denoting by $\bP$ the corresponding p-measured lattice, we set
$Q=\setm{x\in P}{\bve{x}{0}\in S}$, which is a lower subset
of~$P$. The restriction of the p-measure of $\bP$ to $Q\times Q$
is as required.
\end{proof}

The following easy result shows that distributivity cannot be removed from the assumptions of Corollary~\ref{C:PosetRepr}.

\begin{proposition}\label{P:CondImplDistr}
Let $S$ be a \jzs, let $P$ be a poset, and let $\bvi{{}_{-}}{{}_{-}}$ be a $S$-valued p-measure on $P$ satisfying condition \textup{(ii)} of Corollary~\textup{\ref{C:PosetRepr}} such that the subset $\Sigma=\setm{\bve{y}{x}}{x\leq y\text{ in }P}$ join-generates~$S$. Then $S$ is distributive.
\end{proposition}

\begin{proof}
Let $\ba_0,\ba_1,\bb\in S$ such that $\bb\leq\ba_0\vee\ba_1$, we find $\bb_i\leq\ba_i$, for $i<2$, such that $\bb=\bb_0\vee\bb_1$. Suppose first that $\bb\in\Sigma$,
so $\bb=\bve{y}{x}$, for some $x\leq y$ in $P$. By assumption, there are a positive integer
$m$ and a decomposition $x=z_0\leq z_1\leq\cdots\leq z_m=y$ such that for each $i<m$, there exists $\eps(i)\in\set{0,1}$ with $\bve{z_{i+1}}{z_i}\leq\ba_{\eps(i)}$.
Put $\bb_j=\bigvee\famm{\bve{z_{i+1}}{z_i}}{i\in\eps^{-1}\set{j}}$, for all $j<2$. Then $\bb_j\leq\ba_j$ and $\bb=\bve{y}{x}=\bb_0\vee\bb_1$.

In the general case, $\bb=\bigvee\famm{\bc_j}{j<n}$ for a positive integer $n$ and elements
$\bc_0,\dots,\bc_{n-1}\in\Sigma$. By the above paragraph, there are decompositions $\bc_j=\bc_{j,0}\vee\bc_{j,1}$ with $\bc_{j,k}\leq\ba_k$ for all $j<n$ and $k<2$. The elements $\bb_k=\bigvee\famm{\bc_{j,k}}{j<n}$, for $k<2$, are as required.
\end{proof}

The following example shows that the conditions (DB1) and (DB2)
cannot be removed from the assumptions of
Theorem~\ref{T:ConstrMeasDiagr}. The construction is inspired by
the one of the cube $\cD_{\mathrm{c}}$ presented in
\cite[Section~3]{TuWe1}.

\begin{examplepf}\label{Ex:PosetVersDc}
Put $\Lambda=\Pow(3)$ \pup{the three-dimensional cube}
and $\Lambda^*=\Lambda\setminus\set{3}$.
There are a $\Lambda$-indexed diagram
$\cB=\seqm{B_p}{p\in\Lambda}$ of finite Boolean lattices and
\jzue s, whose restriction to~$\Lambda^*$ we denote by $\cB^*$,
and a $\cB^*$-valued normal interval diagram
$\seqm{\bQ_p}{p\in\Lambda^*}$ of finite p-measured lattices that
cannot be extended to any $\cB$-valued normal diagram of
p-measured posets.
\end{examplepf}

\begin{proof}
We first put $B_{\set{0,1,2}}=\Pow(5)$
(where, as usual, $5=\set{0,1,2,3,4}$). Further, we define
elements $\bc_{i,j}$ of $\Pow(5)$, for $i<3$ and $j<4$, by
 \begin{align*}
 \bc_{0,0}&=\set{0,4},&\bc_{0,1}&=\set{3},&
 \bc_{0,2}&=\set{2},&\bc_{0,3}&=\set{1,4};\\
 \bc_{1,0}&=\set{0,4},&\bc_{1,1}&=\set{1,4},&
 \bc_{1,2}&=\set{2},&\bc_{1,3}&=\set{3,4};\\
 \bc_{2,0}&=\set{0,4},&\bc_{2,1}&=\set{1},&
 \bc_{2,2}&=\set{3},&\bc_{2,3}&=\set{2,4}.
 \end{align*}
Observe that the equality $5=\bigcup\famm{\bc_{i,j}}{j<4}$
holds, for all $i<3$.

We shall now define certain subsemilattices of
$\seq{\Pow(5),\cup,\es}$. For $\set{i,j,k}=3$, we define
$B_{\set{i,j}}$ as the \jz-subsemilattice of
$\seq{\Pow(5),\cup,\es}$ generated by the subset
$\set{\bc_{k,0},\bc_{k,1},\bc_{k,2},\bc_{k,3}}$.

Further, for all $i<3$, let $B_{\set{i}}$ be the
\jz-subsemilattice of $\Pow(5)$ generated by
$\set{\ba_i,\bb_i}$, where we put
 \begin{align*}
 \ba_0&=\set{0,1,4},&\bb_0&=\set{2,3,4};\\
 \ba_1&=\set{0,3,4},&\bb_1&=\set{1,2,4};\\
 \ba_2&=\set{0,2,4},&\bb_2&=\set{1,3,4}.
 \end{align*}
At the bottom of the diagram, we put the two-element
semilattice $B_{\es}=\set{\es,5}$. Observe, in particular, that
$5$ is the largest element of $B_p$ for all $p\subseteq 3$.

It is a matter of routine to verify that $B_p$ is a
\jzu-subsemilattice of $B_q$ if $p\subseteq q$, for
all $p,q\subseteq 3$. In that case, we denote by $\varphi_{p,q}$
the inclusion map from~$B_p$ into~$B_q$. Set
 \begin{align*}
 \cB&=\seqm{\seq{B_p,\varphi_{p,q}}}
 {p\subseteq q\text{ in }\Pow(3)},\\
 \cB^*&=\seqm{\seq{B_p,\varphi_{p,q}}}
 {p\subseteq q\text{ in }\Pow(3)\setminus\set{3}}.
 \end{align*}
Let $Q_p$, for $p\in\Lambda^*$, and $P$ be the
lattices diagrammed on Figure~\ref{Fig:DcPieces}.
We observe that $\seqm{Q_p}{p\in\Lambda^*}$ is a normal interval
diagram of finite lattices.
\begin{figure}[htb]
\includegraphics{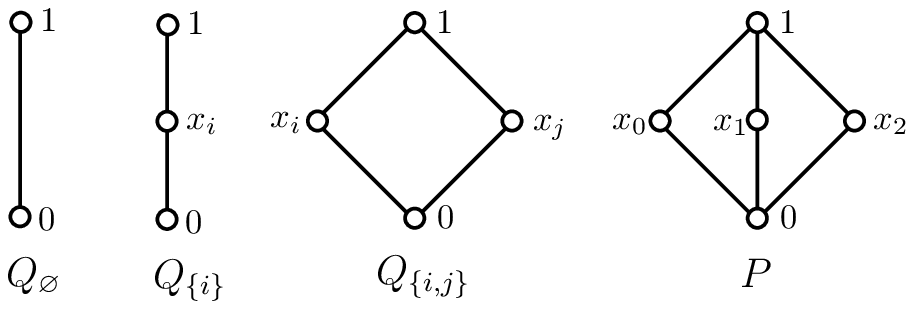}
\caption{The posets $Q_{\es}$, $Q_{\set{i}}$, $Q_{\set{i,j}}$,
and $P$.}
\label{Fig:DcPieces}
\end{figure}
We endow~$Q_{\es}$ with the unique p-measure
$\bvi{{}_{-}}{{}_{-}}_{\es}$ that satisfies
$\bve{1}{0}_{\es}=5$, the largest element of $B_{\es}$. For
$i<3$, we endow $Q_{\set{i}}$ with the unique p-measure
$\bvi{{}_{-}}{{}_{-}}_{\set{i}}$ that satisfies
$\bve{x_i}{0}_{\set{i}}=\ba_i$ and
$\bve{1}{x_i}_{\set{i}}=\bb_i$. Finally, for
$\set{i,j,k}=3$, it is not hard to verify that there exists a
unique p-measure $\bvi{{}_{-}}{{}_{-}}_{\set{i,j}}$ on
$Q_{\set{i,j}}$ such that
$\bvi{x_i}{x_j}_{\set{i,j}}=\bc_{k,1}$ and
$\bvi{x_j}{x_i}_{\set{i,j}}=\bc_{k,2}$.

Suppose that the $\cB^*$-valued diagram
$\seqm{\bQ_p}{p\in\Lambda^*}$ extends to some $\cB$-valued
diagram $\seqm{\bQ_p}{p\in\Lambda}$. Evaluating the Boolean
values in $\bQ_{\set{0,1,2}}$, we obtain
 \begin{align*}
 \bvi{x_0}{x_1}&=\bvi{x_0}{x_1}_{\set{0,1}}=\bc_{2,1}\,,\\
 \bvi{x_1}{x_2}&=\bvi{x_1}{x_2}_{\set{1,2}}=\bc_{0,1}\,,\\
 \bvi{x_0}{x_2}&=\bvi{x_0}{x_2}_{\set{0,2}}=\bc_{1,1}\,,
 \end{align*}
hence, by the triangular inequality,
$\bc_{1,1}\subseteq\bc_{0,1}\cup\bc_{2,1}$, \contr.
\end{proof}

\section{Concluding remarks}
\label{S:Dist}

\subsection{Relation with the V-distances of \cite{RTW}}\label{Su:Vdist}
The main result of the present paper, Theorem~\ref{T:PosetRepr},
is formally similar to \cite[Theorem~7.1]{RTW},
which states that \emph{every distributive \jzs\ is, functorially,
the range of a V-distance of type~$2$ on some set}. By
definition, for a \jzs\ $S$, a $S$-valued distance on a set
$X$ is a map $\delta\colon X\times X\to S$ such that
$\delta(x,x)=0$, $\delta(x,y)=\delta(y,x)$, and
$\delta(x,z)\leq\delta(x,y)\vee\delta(y,z)$, for all
$x,y,z\in X$. Furthermore, $\delta$ satisfies the V-condition of
type~$2$, if for all $\ba,\bb\in S$ and all $x,y\in X$, if
$\delta(x,y)=\ba\vee\bb$, then there are $u,v\in X$ such that
$\delta(x,u)\vee\delta(v,y)\leq\ba$ and $\delta(u,v)\leq\bb$. (The `V-condition' is named so after Hans Dobbertin's work in~\cite{Dobb83}.) As
every distance on a set $X$ is obviously a p-measure on $X$
viewed as a discrete poset, the problem of functorially lifting
distributive \jzs s by p-measures does not appear as
difficult. The main problems encountered in the present work were
(1) to get our posets \emph{connected} (which is the case here as they
are \ms s), and (2) to get the subset $\setm{\bve{y}{x}}{x\leq y\text{ in }P}$
join-generating the semilattice $S$ under consideration.

\subsection{Representation of distributive semilattices by majority algebras}\label{Su:MajAlg}
To the author's knowledge, Corollary~\ref{C:PosetRepr} is, so far, the only existing representation result that is \emph{specific} to distributive \jzs s. Unlike the Gr\"atzer-Schmidt Theorem, it is not a lifting result of \jzs s with respect to the~$\Conc$ functor---the functor under consideration, namely~$\Pi$ (cf. Subsection~\ref{Su:Lifting}), is more complicated to describe. One remaining hope after the negative result of~\cite{CLP} is whether every distributive \jzs\ is isomorphic to $\Conc A$ for some algebra~$A$ generating a \emph{congruence-distributive} variety (cf. \cite[Problem~2]{CLP}). For instance, \emph{is every distributive \jzs\ isomorphic to $\Conc M$, for some majority algebra~$M$}? (A \emph{majority algebra} is a nonempty set endowed with a ternary operation~$m$ that satisfies the identities
$m(\sx,\sx,\sy)=m(\sx,\sy,\sx)=m(\sy,\sx,\sx)=\sx$.) Our hope is that the poset-theoretical methods used in the present paper could provide a stepping stone towards such a result.

\subsection{Lifting finite diagrams of finite distributive \jzs s}\label{Su:FinDiagr}
It is still an open problem whether every diagram~$\vec{D}$ of finite \jzs s and \jzh s, indexed by a finite lattice, can be lifted, with respect to the~$\Conc$ functor, by a diagram of (finite?) lattices (cf. \cite[Problem~4]{CLP}). Applying Theorem~\ref{T:PosetRepr} to the diagram of \jzus s obtained by adding a largest element to each object in~$\vec{D}$ and extending the transition maps accordingly, and then restricting the posets as at the end of the proof of Corollary~\ref{C:PosetRepr}, gives the weaker result that~$\vec{D}$ can be lifted, with respect to the~$\Pi$ functor (cf. Subsection~\ref{Su:Lifting}), by a diagram in~$\VPmeas$. The posets thus obtained may be thought of as `skeletons' of the lattices that would appear in a (hypothetical) lifting of~$\vec{D}$ with respect to~$\Conc$.

\end{document}